\documentclass[a4paper,11pt]{article}
\usepackage[T1]{fontenc}
\usepackage[latin9]{inputenc}
\usepackage{palatino}
\usepackage{hyperref}
\hypersetup{colorlinks=true, linkcolor=blue,  anchorcolor=blue, citecolor=blue, filecolor=blue, menucolor=blue, urlcolor=blue}
\usepackage[font={footnotesize,sf},labelfont=bf]{caption}
\usepackage{amsmath,amssymb,amsfonts,mathrsfs,amsthm}
\usepackage{graphicx,pdfpages,enumitem,cite}
\usepackage{algorithm,algorithmic}
\theoremstyle{plain}

\theoremstyle{plain} 
 
\theoremstyle{definition} 
 
\usepackage{geometry}
\title{Data-Driven Modeling and Prediction of Non-Linearizable Dynamics via
Spectral Submanifolds}
\author{Mattia Cenedese$^{1}$, Joar Axås$^{1}$, Bastian Bäuerlein$^{2,3}$, \\
Kerstin Avila$^{2,3}$ and George Haller$^{1,}$\thanks{Corresponding author email: georgehaller@ethz.ch}\vspace{3.5mm}\\
$^{1}$Institute for Mechanical Systems, ETH Zürich, \\Leonhardstrasse
21, 8092 Zürich, Switzerland\vspace{1.5mm}\\
$^{2}$University of Bremen, Faculty of Production Engineering, \\Badgasteiner
Strasse 1, 28359, Bremen, Germany\vspace{1.5mm}\\
 $^{3}$Leibniz Institute for Materials Engineering IWT, \\Badgasteiner
Strasse 3, 28359, Bremen, Germany
}
\date{}
\begin{document}
\maketitle
\begin{abstract}
We develop a methodology to construct low-dimensional predictive models
from data sets representing essentially nonlinear (or \emph{non-linearizable})
dynamical systems with a hyperbolic linear part that are subject to
external forcing with finitely many frequencies. Our data-driven,
sparse, nonlinear models are obtained as extended normal forms of
the reduced dynamics on low-dimensional, attracting spectral submanifolds
(SSMs) of the dynamical system. We illustrate the power of data-driven
SSM reduction on high-dimensional numerical data sets and experimental
measurements involving beam oscillations, vortex shedding and sloshing
in a water tank. We find that SSM reduction trained on unforced data
also predicts nonlinear response accurately under additional external
forcing.
\end{abstract}

\section{Introduction\label{sec:Introduction}}

Low-dimensional reduced models of high-dimensional nonlinear dynamical
systems are critically needed in various branches of applied science
and engineering. Such simplified models would significantly reduce
computational costs and enable physical interpretability, design optimization
and efficient controllability. As of yet, however, no generally applicable
procedure has emerged for the reliable and robust identification of
nonlinear reduced models.

Instead, the most broadly used approach to reducing nonlinear dynamical
systems has been a fundamentally linear technique, the proper orthogonal
decomposition (POD), followed by a Galerkin projection \cite{holmes12,awrejcewicz04,lu19}.
Projecting the full dynamics to the most energetic linear modes, POD
requires the knowledge of the governing equations of the system and
hence is inapplicable when only data is available. As purely data-based
alternatives, machine learning methods are broadly considered and
tested in various fields \cite{brunton2016,mohamed18,daniel20,calka21}.
While the black-box approach of machine learning might often seem
preferable to a detailed nonlinear analysis, the resulting neural
network models require extensive tuning, lack physical interpretability,
generally perform poorly outside their training range and tend to
be unnecessarily complex \cite{loiseau20}. This has inspired a number
of approaches that seek a blend of machine learning with a priori
information about the underlying physics \cite{karniadakis21,li2021}.
Still within the realm of machine learning, sparse regression has
also shown promise in approximating the right-hand sides of low-dimensional,
simple dynamical systems with functions taken from a preselected library
\cite{brunton2016}. Another recent approach is cluster-based network
modeling, which uses the toolkit of network science and statistical
physics for modeling nonlinear dynamics \cite{fernex21}.

A popular alternative to POD and machine learning is the dynamic mode
decomposition (DMD) \cite{schmid10}, which approximates directly
the observed system dynamics. The original DMD and its later variants
fit a linear dynamical system to temporally evolving data, possibly
including further functions of the original data, over a given finite
time interval \cite{kutz16}. DMD provides an appealingly simple
yet powerful algorithm to infer a local model near steady states where
the nonlinear dynamics is always approximately linear. This linear
model is also more globally valid if constructed over observables
lying in a span of some eigenfunctions of the Koopman operator, which
maps observables evaluated over initial states into their evaluations
over current states \cite{rowley09,mezic13,mauroy20}. This relationship
between DMD and the Koopman operator has motivated an effort to machine-learn
Koopman eigenfunctions from data in order to linearize nonlinear dynamical
systems globally on the space of their observables \cite{lusch18,otto19,kaiser21}.

Finding physically realizable observables that fall in a Koopman eigenspace
is, however, often described as challenging or difficult \cite{page19}.
A more precise assessment would be that such a find is highly unlikely,
given that the probability of any countable set of a priori selected
observables falling in any Koopman eigenspace is zero. In addition,
those eigenspaces can only be determined explicitly in simple, low-dimensional
systems. In practice, therefore, DMD can only provide a justifiable
model near an attracting fixed point of a dynamical system. While
Koopman modes still have the potential to linearize the observer dynamics
on larger domains, those domains cannot include more than one attracting
or repelling fixed point \cite{page19,brunton16,kaiser21}. Indeed,
DMD and Koopman mode expansions fail to converge outside neighborhoods
of fixed points even in the simplest, one-dimensional nonlinear systems
with two fixed points \cite{bagheri13,page19}. In summary, while
these data-driven model reduction methods are powerful and continue
to inspire ongoing research, their applicability is limited to locally
linearized systems and globally linearizable nonlinear systems, such
as the Burgers equation \cite{page18}. 

The focus of this paper is the development of data-driven, simple
and predictive reduced-order models for essentially nonlinear dynamical
systems, i.e., \emph{non-linearizable systems}. Determining exact
linearizability conclusively from data is beyond reach. In contrast,
establishing that a dynamical system is non-linearizable in a domain
of interest is substantially simpler: one just needs to find an indication
of coexisting isolated stationary states in the data. By an isolated
stationary state, we mean here a compact and connected invariant set
with an open neighborhood that contains no other compact and connected
invariant set. Examples of such stationary states include hyperbolic
fixed points, periodic orbits, invariant spheres and quasiperiodic
tori; closures of homoclinic orbits and heteroclinic cycles; and chaotic
attractors and repellers. If a data set indicates the coexistence
of any two sets from the above list, then the system is conclusively
non-linearizable in the range of the available data. Specifically,
there will be no homeomorphism (continuous transformation with a continuous
inverse) that transforms the orbits of the underlying dynamical system
into those of a linear dynamical system. While this is a priori clear
from dynamical systems theory, several studies have specifically confirmed
a lack of convergence of Koopman-mode expansions already for the simplest
case of two coexisting fixed points, even over subsets of their domain
of attraction or repulsion \cite{bagheri13,page19}.

\begin{figure}[t]
\includegraphics[width=1\textwidth]{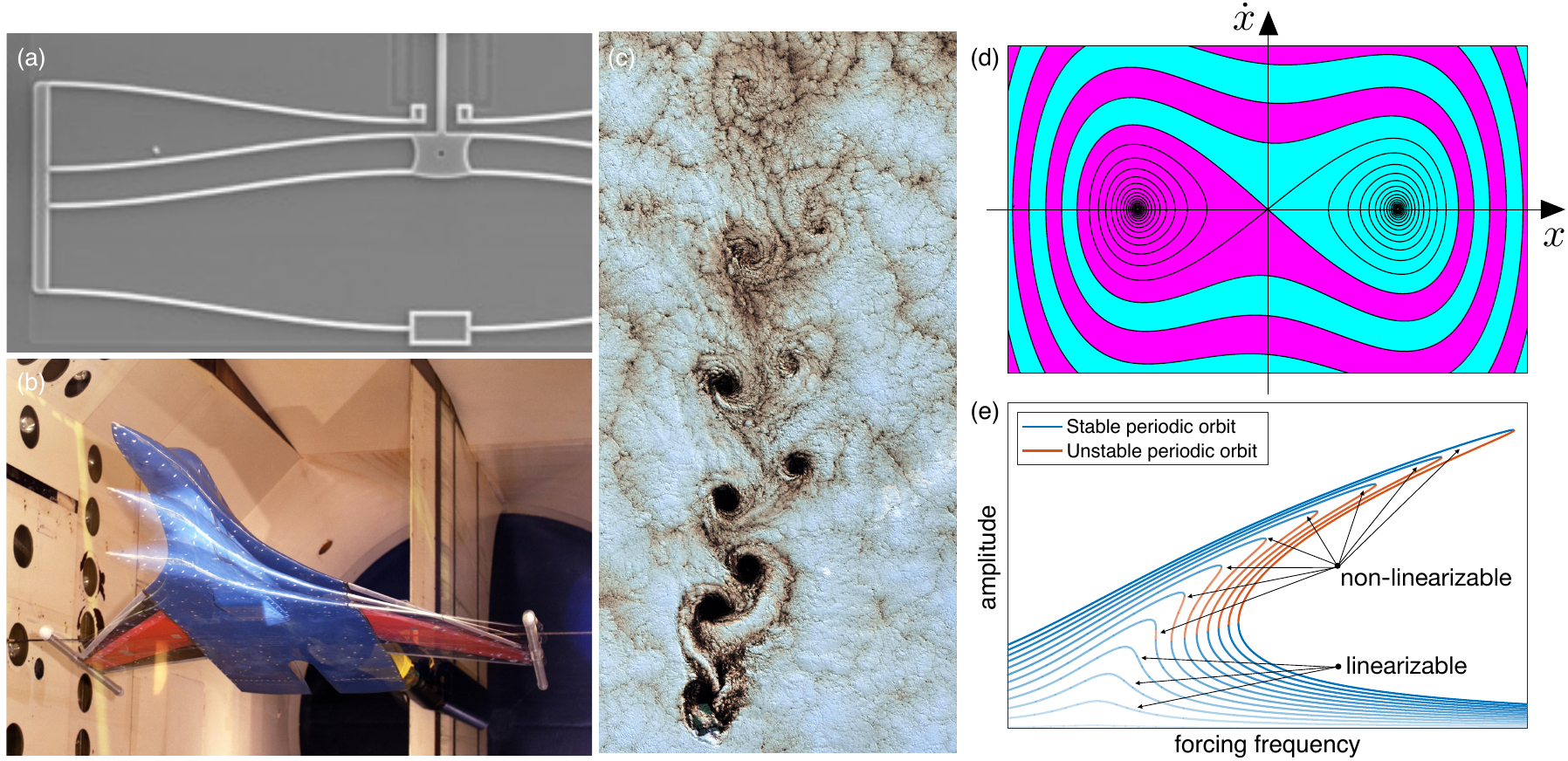}
\caption{Examples of non-linearizable systems: (a) Snap-through instability
of a micro-electro-mechanical (MEMS) device with three coexisting
equilibria (Sandia National Laboratories) (b) Wind-tunnel flutter
of an airplane prototype, involving a fixed point and coexisting limit
cycles (NASA Langley Research Center). (c) Swirling clouds behind
an island in the Pacific ocean, forming a vortex street with coexisting
isolated hyperbolic and elliptic trajectories for the dynamical system
describing fluid particle motion (USGS/NASA). (d) Phase portrait of
the damped, double-well Duffing oscillator $\ddot{x}+\dot{x}-x+\beta x^{3}=0$
with $\beta<0$, the most broadly used model for nonlinear systems
with coexisting domains of attraction (colored), such as the MEMS
device in plot (a). (e) Nonlinear response amplitude ($|x(t)|_{\mathrm{max}}$)
in the forced-damped, single-well Duffing oscillator, $\ddot{x}+\dot{x}+x+\beta x^{3}=\cos\omega t$
with $\beta>0$, under variations of the forcing frequency $\omega$.
Coexisting stable and unstable periodic responses show non-linearizable
dynamics conclusively for this classic model.}
\label{fig:nonlinearizablesystems}
\end{figure}

Non-linearizable systems are ubiquitous in science, technology and
nature. Beyond the well-known examples of chaotic dynamical systems
and turbulent fluid flows \cite{holmes12}, any bifurcation phenomenon,
by definition, involves coexisting steady states and hence is automatically
non-linearizable. Indeed, aerodynamic flutter \cite{dowell70}, buckling
of beams and shells \cite{abramian20}, bistable microelectromechanical
systems \cite{podder2016}, traffic jams \cite{orosz06} or even
tipping points in climate change \cite{ashwin12} are all fundamentally
non-linearizable, just to name a few. Figure \ref{fig:nonlinearizablesystems}
shows some examples of non-linearizable systems emerging in technology,
nature and scientific modeling.

We will show here that a collection of classic and recent mathematical
results from nonlinear dynamical systems theory enables surprisingly
accurate and predictive low-dimensional modeling from data for a number
of non-linearizable phenomena. Our construct relies on the recent
theory of \emph{spectral submanifolds} (SSMs), the smoothest invariant
manifolds that act as nonlinear continuations of non-resonant eigenspaces
from the linearization of a system at a stationary state (fixed point,
periodic orbit or quasiperiodic orbit; \cite{haller16}). Using appropriate
SSM embeddings \cite{whitney44,stark97,stark99} and an extended
form of the classic normal form theory \cite{guckenheimer83}, we
obtain sparse dynamical systems describing the reduced dynamics on
the slowest SSMs of the system, which are normally hyperbolic and
hence robust under perturbations \cite{fenichel71}. 

We construct the extended normal form within the slowest SSM as if
the eigenvalues of the linearized dynamics within the SSM had zero
real parts, although that is not the case. As a result, our normalization
procedure will not render the simplest possible (linear) normal form
for the SSM dynamics, valid only near the underlying isolated stationary
state. Instead, our procedure yields a sparsified nonlinear, polynomial
normal form on a larger domain of the SSM that can also capture nearby
coexisting stationary states. This fully data-driven normalization
algorithm learns the normal form transformation and the coefficients
of the normal form simultaneously by minimizing an appropriately defined
conjugacy error between the unnormalized and normalized SSM dynamics.

For a generic observable of an oscillatory dynamical system without
an internal resonance, a two-dimensional data-driven model calculated
on the slowest SSM of the system turns out to capture the correct
asymptotic dynamics. Such an SSM-reduced model is valid on domains
in which the nonlinearity and any possible external forcing are strong
enough to create non-linearizable dynamics, yet are still moderate
enough to render the eigenspace of the linear system relevant. {\scriptsize{}}
More generally, oscillatory systems with $m$ independent internal
resonances in their spectrum can be described by reduced models on
$2\left(m+1\right)$-dimensional SSMs. In both the resonant and the
non-resonant cases, the models can be refined by increasing the degree
of their nonlinearity rather than by increasing their dimension. As
we show in examples, the resulting SSM-based models are explicit,
deterministic and even have the potential to predict system behavior
outside the range of the training data away from bifurcations. Most
importantly, we find that the models also accurately predict forced
response, even though they are only trained on data collected from
unforced systems. 

We illustrate the power of data-driven SSM-reduced models on high-dimensional
numerically generated data sets and on experimental data. These and
further examples are also available as \textsc{MATLAB}\textsuperscript{®}
live scripts, which are part of a general open-source package, \texttt{SSMLearn,}
that performs this type of model reduction and prediction for arbitrary
data sets.

\section{Results}

\subsection{Spectral submanifolds and their reduced dynamics}

A recent result in dynamical systems is that all eigenspaces (or \emph{spectral
subspaces}) of linearized systems admit unique nonlinear continuations
under well-defined mathematical conditions. Specifically, \emph{spectral
submanifolds} (SSMs), as defined by \cite{haller16}, are the unique
smoothest invariant manifolds that serve as nonlinear extensions of
spectral subspaces under the addition of nonlinearities to a linear
system. The SSM formulation and terminology we use here is due to
\cite{haller16}; the Methods section \ref{subsec:Existence-of-SSMs}
discusses the history of these results and further technical details. 

We consider $n$-dimensional dynamical systems of the form 
\begin{equation}
\dot{\mathbf{x}}=\mathbf{A}\mathbf{x}+\mathbf{f}_{0}(\mathbf{x})+\epsilon\mathbf{f}_{1}(\mathbf{x},\boldsymbol{\Omega}t;\epsilon),\qquad\mathbf{f}_{0}(\mathbf{x})=\mathcal{O}(\left|\mathbf{x}\right|^{2}),\qquad0\leq\epsilon\ll1,\label{eq:1storder_system}
\end{equation}
with a constant matrix $\mathbf{A}\in\mathbb{R}^{n\times n},$ and
with class $C^{r}$ functions $\mathbf{f}_{0}\colon\mathcal{U}\to\mathbb{R}^{n}$
and $\mathbf{f}_{1}\colon\mathcal{U}\times\mathbb{T}^{\ell}\to\mathbb{R}^{n}$,
where $\mathbb{T}^{\ell}=S^{1}\times\ldots\times S^{1}$ is the $\ell$-dimensional
torus. The elements of the frequency vector $\boldsymbol{\Omega}\mathbb{\in}\mathbb{R}^{\ell}$
are rationally independent, and hence the function $\mathbf{f}_{1}$
is quasiperiodic in time. The assumed degree of smoothness for the
right-hand side of (\ref{eq:1storder_system}) is $r\in\mathbb{N}^{+}\cup\left\{ \infty,a\right\} $,
with $a$ referring to analytic. The small parameter $\epsilon$ signals
that the forcing in system (\ref{eq:1storder_system}) is moderate
so that the structure of the autonomous part is still relevant for
the full system dynamics. Rigorous mathematical results on SSMs are
proven for small enough $\epsilon$, but continue to hold in practice
for larger values of $\epsilon$ as well, as we will see in examples.
Note that eq. \eqref{eq:1storder_system} describes equations of motions
of physical oscillatory systems. It does not cover phenomenological
models of phase oscillators, such as the Kuramoto model \cite{kuramoto84}.

The eigenvalues $\lambda_{j}=\alpha_{j}+\mathrm{i}\omega_{j}\in\mathbb{C}$
of\textbf{ $\mathbf{A}$}, with multiplicities counted, are ordered
based on their real parts, $\mathrm{Re}\lambda_{j}$, as
\begin{equation}
\mathrm{Re}\lambda_{n}\leq\mathrm{Re}\lambda_{n-1}\leq\ldots\ldots\leq\mathrm{Re}\lambda_{1}.\label{eq:eigenvalue_ordering}
\end{equation}
Their corresponding real \emph{modal subspaces} (or eigenspaces),
$E_{j}\subset\mathrm{\mathbb{R}}^{n}$ , are spanned by the imaginary
and real parts of the corresponding eigenvectors and generalized eigenvectors
of $\mathbf{A}$. To analyze typical systems, we assume that $\mathrm{Re}\lambda_{j}=\alpha_{j}\neq0$
holds for all eigenvalues, i.e., $\mathbf{x}=\mathbf{0}$ is a hyperbolic
fixed point for $\epsilon=0$. 

\begin{figure}[t]
\includegraphics[width=1\textwidth]{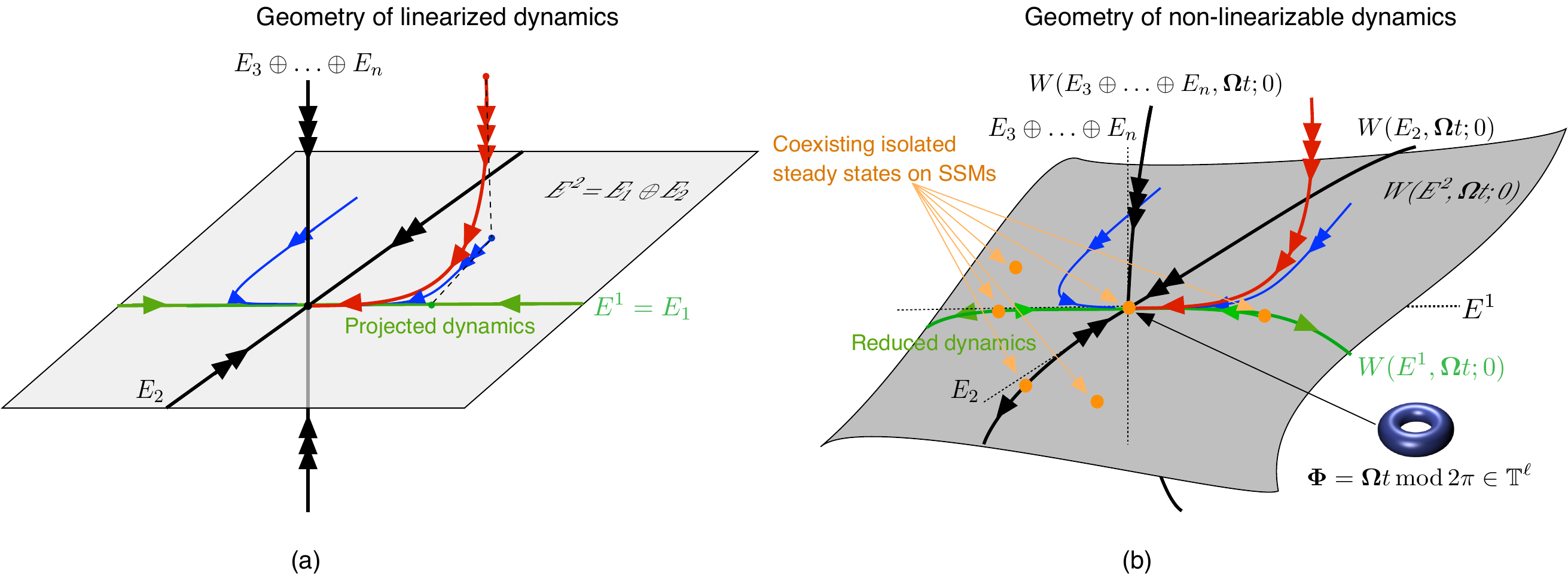}
\caption{(a) Reduction of linear dynamics via Galerkin projection. The slowest
spectral subspace, $E^{1}=E_{1}$ (green), and the modal subspace,
$E_{2}$ (black), span together the second slowest spectral subspace,
$E^{2}=E_{1}\oplus E_{2}$ . The full dynamics (red curve) can be
projected onto $E^{1}$ to yield a reduced slow model without transients.
Projection of the full dynamics onto $E^{2}$ (blue curve) yields
a reduced model that also captures the slowest decaying transient.
Further, faster-decaying transients can be captured by projections
onto slow spectral subspaces, $E^{k}$, with $k>1$. (b) Reduction
of non-linearizable dynamics via restriction to spectral submanifolds
(SSMs) in the $\epsilon=0$ limit of nonlinear, non-autonomous systems
forced with $\ell$ frequencies. An SSM, $W(E,\boldsymbol{\Omega}t;0)$,
is the unique, smoothest, nonlinear continuation of a non-resonant
spectral subspace $E$. Specifically, the slowest SSM, $W(E^{k},\boldsymbol{\Omega}t;0)$
(green), is the unique, smoothest, nonlinear continuation of the slowest
spectral subspace, $E^{k}$. Non-linearizability of the full dynamics
follows if isolated stationary states coexist on al least one of the
SSMs. The time-quasiperiodic SSMs for $\epsilon>0$, denoted $W(E,\boldsymbol{\Omega}t;\epsilon)$,
are not shown here but they are $\mathcal{O}(\epsilon)$ $C^{r}$-close
to the structures shown, as discussed by \cite{haller16}.}\label{fig:linear vs nonlinear dynamics}
\end{figure}

A \emph{spectral subspace} $E_{j_{1},\ldots,j_{q}}$ is a direct sum
\emph{
\begin{equation}
E_{j_{1},\ldots,j_{q}}=E_{j_{1}}\oplus E_{j_{2}}\oplus\ldots\oplus E_{j_{q}}\label{eq:spectral subspace}
\end{equation}
}of an arbitrary collection of modal subspaces, which is always an
invariant subspace for the linear part of the dynamics in (\ref{eq:1storder_system}).
Classic examples of spectral subspaces are the stable and unstable
subspaces, comprising all modal subspaces with $\mathrm{Re}\lambda_{k}<0$
and $\mathrm{Re}\lambda_{k}>0$, respectively. Projections of the
linearized system onto the nested hierarchy of \emph{slow spectral
subspaces}, 
\begin{equation}
E^{1}\subset E^{2}\subset E^{3}\subset\ldots,\qquad E^{k}:=E_{1,\ldots,k},\quad k=1,\ldots,n,\label{eq:hierarchy of slow spectral subspaces}
\end{equation}
provide exact reduced-order models for the linearized dynamics over
an increasing number of time scales under increasing $k$, as sketched
in Fig. \ref{fig:linear vs nonlinear dynamics}a. This is why a Galerkin
projection onto $E^{k}$ is an exact model reduction procedure for
linear systems, whose accuracy can be increased by increasing $k$.
A fundamental question is whether nonlinear analogues of spectral
subspaces continue to organize the dynamics under the addition of
nonlinear and time-dependent terms in the full system (\ref{eq:1storder_system}).

Let us fix a specific spectral subspace $E=E_{j_{1},\ldots,j_{q}}$
within either the stable or the unstable subspace. If $E$ is non-resonant
(i.e., no nonnegative, low-order, integer linear combination of the
spectrum of $\mathbf{A}\vert_{E}$ is contained in the spectrum of
$\mathbf{A}$ outside $E$), then $E$ has infinitely many nonlinear
continuations in system (\ref{eq:1storder_system}) for $\epsilon$
small enough \cite{haller16}. These invariant manifolds are of smoothness
class $C^{\Sigma(E)}$, with the spectral quotient $\Sigma(E)$ measuring
the ration of the fastest decay exponent outside $E$ to the slowest
decay exponent inside $E$ (see eq. (\eqref{eq:absolute spectral quotient})
of the Methods section \ref{subsec:Existence-of-SSMs}). All
such manifolds are tangent to $E$ for $\epsilon=0$, have the same
quasiperiodic time dependence as $\mathbf{f}_{1}$ does and have a
dimension equal to that of $E$. 

Of these infinitely may invariant manifolds, however, there will be
a unique smoothest one, the \emph{spectral submanifold} (SSM) of $E$,
denoted $W(E,\boldsymbol{\Omega}t;\epsilon)$. This manifold is $C^{r}$
smooth if $r>\Sigma(E)$ and can therefore be approximated more accurately
than the other infinitely many nonlinear continuations of $E$. In
particular, SSMs have convergent Taylor expansions if the dynamical
system \eqref{eq:1storder_system} is analytic ($r=a$). Then the
reduced dynamics on a slow SSM, $E^{k}$, can be approximated with
arbitrarily high accuracy using arbitrarily high-order Taylor expansions,
without ever increasing the dimension of $E^{k}$ (see Fig. \ref{fig:linear vs nonlinear dynamics}b).
Such an approximation for dynamical systems with known governing equations
is now available for any required order of accuracy via the open source
\textsc{MATLAB}\textsuperscript{®} package \texttt{SSMTool} \cite{jain2021}.
In contrast, reduced models obtained from projection-based procedures
can only be improved by increasing their dimensions. 

The nearby coexisting stationary states in Fig. \ref{fig:linear vs nonlinear dynamics}
happen to be contained in the SSM. In specific examples, however,
these states may also be off the SSM, contained instead in one of
the infinitely many additional nonlinear continuations, $\tilde{W}(E,\boldsymbol{\Omega}t;\epsilon)$,
of the spectral subspace $E$. The Taylor expansion of the dynamics
on $\tilde{W}(E,\boldsymbol{\Omega}t;\epsilon)$ and $W(E,\boldsymbol{\Omega}t;\epsilon)$
are, however, identical up to order $\Sigma(E)$. Therefore, the reduced
models we will compute on the SSM $W(E,\boldsymbol{\Omega}t;\epsilon)$
also correctly capture the nearby stationary states on $\tilde{W}(E,\boldsymbol{\Omega}t;\epsilon)$,
as long as the polynomial order of the model stays below $\Sigma(E)$.
In large physical systems, this represents no limitation, given that
$\Sigma(E)\gg1$.

\subsection{Embedding SSMs via generic observables\label{subsec:Embedding-of-SSMs}}

If at least some of the real parts of the eigenvalues in \eqref{eq:eigenvalue_ordering}
are negative, then longer-term trajectory data for system (\ref{eq:1storder_system})
will be close to an attracting SSM, as illustrated in Fig. \ref{fig:linear vs nonlinear dynamics}b.
This is certainly the case for data from experiments that are run
until a nontrivial, attracting steady state emerges (see, e.g., Fig.
\ref{fig:nonlinearizablesystems}e). Measurements of trajectories
in the full phase space, however, are seldom available from such experiments.
Hence, if data about system (\ref{eq:1storder_system}) is only available
from observables, the construction of SSMs and their reduced dynamics
has to be carried out in the space of those observables.

An extended version of Whitney's embedding theorem guarantees that
almost all (in the sense of prevalence) smooth observable vectors
$\mathbf{y}(\mathbf{x})=(y_{1}(\mathbf{x}),,...,y_{p}(\mathbf{x}))\in\mathbb{R}^{p}$
provide an embedding of a compact subset $\mathcal{C}\subset W(E,\boldsymbol{\Omega}t;\epsilon)$
of a $d$-dimensional SSM, $W(E,\boldsymbol{\Omega}t;\epsilon)$,
into the observable space $\mathbb{R}^{p}$ for high enough $p$.
Specifically, if we have $p>2(d+\ell)$ simultaneous and independent
continuous measurements, $\mathbf{y}(\mathbf{x})$, of the $p$ observables,
then almost all maps $\mathbf{y}\colon\mathcal{C}\to\mathbb{R}^{p}$
are embeddings of $\text{\ensuremath{\mathcal{C}} }$ \cite{sauer1991},
and hence the top right plot of Fig. \ref{fig:ssmlearn} is applicable
with probability one.

In practice, we may not have access to $p>2(d+\ell)$ independent
observables and hence cannot invoke Whitney's theorem. In that case,
we invoke the Takens delay embedding theorem \cite{takens81}, which
covers observable vectors built from $p$ uniformly sampled, consecutive
measured instances of a single observable. More precisely, if $s(t)$
is a generic scalar quantity measured at times $\Delta t$ apart,
then the observable vector for delay-embedding is formed as $\mathbf{y}(t)=\left(s(t),s(t+\Delta t),...,s\left(t+(p-1)\Delta t\right)\right)\in\mathbb{R}^{p}$.
We discuss the embedding, $\mathcal{M}_{0}\subset\mathbb{R}^{p}$,
of an autonomous SSM, $W(E,\boldsymbol{\Omega}t_{0};0)$, in the observable
space $\mathbb{R}^{p}$ in more detail in the Methods section
\ref{subsec:SSM-geometry-in-observable-space}.

\subsection{Data-driven extended normal forms on SSMs}

Once the embedded SSM, $\mathcal{M}_{0}$, is identified in the observable
space, we seek to learn the reduced dynamics on $\mathcal{M}_{0}$.
An emerging requirement for learning nonlinear models from data has
been model sparsity \cite{brunton2016}, without which the learning
process would be highly sensitive. The dynamics on $\mathcal{M}_{0}$,
however, is inherently non-sparse, which suggests that we learn its
Poincaré normal form \cite{poincare1892} instead. This classic normal
form is the simplest polynomial form to which the dynamics can be
brought via successive, near-identity polynomial transformations of
increasing order. 

Near the origin on a slow SSM, however, this simplest polynomial form
is just the restriction of the linear part of system \eqref{eq:1storder_system}
to $\mathcal{M}_{0}$, as long as infinitely many non-resonance conditions
are satisfied for the operator $\mathbf{A}$ \cite{sternberg1958}.
The Poincaré normal form on $\mathcal{M}_{0}$ would, therefore, only
capture the low-amplitude, linearized part of the slow SSM dynamics. 

To construct an SSM-reduced model for non-linearizable dynamics, we
use \emph{extended normal forms}. This idea is motivated by normal
forms used in the study of bifurcations of equilibria on center manifolds
depending on parameters \cite{guckenheimer83,murdock2003}. In that
setting, the normal form transformation is constructed at the bifurcation
point where the system is non-linearizable by definition. The same
transformation is then used away from bifurcations, even though the
normal form of the system would be linear there. One, therefore, gives
up the maximal possible simplicity of the normal form but gains a
larger domain on which the normal form transformation is invertible
and hence captures truly nonlinear dynamics. In our setting, there
is no bifurcation at $\mathbf{x}=\mathbf{0}$, but we nevertheless
construct our normal form transformation as if the eigenvalues corresponding
to the slow subspace $E$ were purely imaginary. This procedure leaves
additional, near-resonant terms in the SSM-reduced normal form, enhancing
the domain on which the transformation is invertible and hence the
normal form is valid.

We determine the normal form coefficients directly from data via the
minimization of a conjugacy error (see the Methods section).
This least-square minimization procedure renders simultaneously the
best fitting normal form coefficients and the best fitting normal
form transformation. As we will find in a specific example, this data-driven
procedure can yield accurate reduced models even beyond the formal
domain of convergence of equation-driven normal forms. 

The simplest extended normal form on a slow SSM of an oscillatory
system arises when the underlying spectral subspace $E$ corresponds
to a pair of complex conjugate eigenvalues. Writing in polar coordinates
and truncating at cubic order, \cite{ponsioen2018} finds this normal
form on the corresponding two-dimensional, autonomous SSM, $\mathcal{M}_{0}$,
to be
\begin{equation}
\begin{aligned}\dot{\rho} & =\alpha_{0}\rho+\beta\rho^{3},\\
\dot{\theta} & =\omega_{0}+\gamma\rho^{2}.
\end{aligned}
\label{eq:cubicnormalform}
\end{equation}
This equation is also known as the Stuart--Landau equation arising
in the unfolding of a Hopf bifurcation \cite{landau1944,stuart1960,fujimura1997}. 

The dynamics of (\ref{eq:cubicnormalform}) is characteristically
non-linearizable when $\alpha_{0}\beta<0$, given that a limit cycle
coexists with the $\rho=0$ fixed point in that case. Further coexisting
steady states will arise when forcing is added to the system, as we
discuss in the next section. We note that the cubic normal form on
two-dimensional SSMs has also been approximated from data in \cite{szalai02}.
That non-sparse procedure fits the full observer dynamics to a low-dimensional,
discrete polynomial dynamical system, then performs an analytic SSM
reduction and a classic normal form transformation on the SSM. 

For higher accuracy, the extended normal form on an oscillatory SSM
of dimension $2m$ is of the form
\begin{equation}
\begin{aligned}\dot{\rho}_{j} & =\alpha_{j}(\boldsymbol{\rho},\boldsymbol{\theta})\rho_{j},\\
\dot{\theta}_{j} & =\omega_{j}(\boldsymbol{\rho},\boldsymbol{\theta}),
\end{aligned}
\,\,\,\,\,\,\,j=1,2,...m,\,\,\,\,\,\,\boldsymbol{\rho}\in\mathbb{R}_{+}^{m},\,\,\,\,\,\,\boldsymbol{\theta}\in\mathbb{T}^{m}.\label{eq:oscillatorsnormalform}
\end{equation}
If the linearized frequencies are non-resonant, then the functions
$\alpha_{j}$ and $\omega_{j}$ only depend on $\boldsymbol{\rho}$
\cite{ponsioen2018}. Our numerical procedure determines these functions
up to the necessary order that ensures a required accuracy for the
reduced-order model on the SSM. This is illustrated schematically
for a four-dimensional slow SSM ($m=2)$ in the bottom right plot
of Fig. \ref{fig:ssmlearn}.

\subsection{Predicting forced dynamics from unforced data\label{subsec:Prediction-of-forced-response}}

With the normalized reduced dynamics \eqref{eq:oscillatorsnormalform}
on the embedded SSM, $\mathcal{M}_{0}$, at hand, we can also make
predictions for the dynamics of the embedded quasiperiodic SSM, $\mathcal{M}_{\epsilon}(\boldsymbol{\Omega}t)$,
of the full system \eqref{eq:1storder_system}. This forced SSM is
guaranteed to be an $\mathcal{O}(\epsilon)$ $C^{r}$-close perturbation
of $\mathcal{M}_{0}$ for moderate external forcing amplitudes. A
strict proof of this fact is available for small enough $\epsilon>0$
\cite{haller16}, but as our examples will illustrate, the smooth
persistence of the SSM, $\mathcal{M}_{\epsilon}(\boldsymbol{\Omega}t)$,
generally holds for all moderate $\epsilon$ values in practice. Such
moderate forcing is highly relevant in a number of technological settings,
including system identification in structural dynamics and fluid-structure
interactions, where the forcing must be moderate to preserve the integrity
of the structure.

We discuss the general extended normal form on $\mathcal{M}_{\epsilon}(\boldsymbol{\Omega}t)$
in the Methods section \ref{subsec:SSM-dynamics-via-extended normal forms}.
In the simplest and most frequent special case, the external forcing
is periodic ($\ell=1$) and $\mathcal{M}_{\epsilon}(\Omega t)$ is
the embedding of the slowest, two-dimensional SSM corresponding to
a pair of complex conjugate eigenvalues. Using the modal forcing amplitude
$f_{1,1}$ and modal phase shift $\phi_{1,1}$ in the general normal
form \eqref{eq:quasiperiodicSSMforcing}, \cite{breunung2018} introduces
the new phase coordinate $\psi=\theta-\Omega t-\phi_{1,1}$ and lets
$f=f_{1,1}$, $\alpha=\alpha_{1}$, $\omega=\omega_{1}$ to obtain
the planar, autonomous, extended normal form on $\mathcal{M}_{\epsilon}(\Omega t)$
as
\begin{equation}
\begin{aligned}\dot{\rho} & =\alpha(\rho)\rho+f\sin\psi,\\
\dot{\psi} & =\omega(\rho)-\Omega+\frac{f}{\rho}\cos\psi
\end{aligned}
\label{eq:periodicSSMforcing}
\end{equation}
at leading order in $\epsilon$. All stable and unstable periodic
responses on the SSM are fixed points of system \eqref{eq:periodicSSMforcing},
with their amplitudes $\rho_{0}$ and phases $\psi_{0}$ satisfying
the equations 
\begin{equation}
\Omega=\omega(\rho_{0})\pm\sqrt{\frac{f^{2}}{\rho_{0}^{2}}-\alpha^{2}(\rho_{0})},\quad\psi_{0}=-\tan^{-1}\left[\frac{\alpha\left(\rho_{0}\right)}{\omega\left(\rho_{0}\right)-\Omega}\right].\label{eq:closedformFRC}
\end{equation}

The first analytic formula in \eqref{eq:closedformFRC} predicts the
\emph{forced response curve} (FRC) of system \eqref{eq:1storder_system},
i.e., the relationship between response amplitude, forcing amplitude
and forcing frequency, from the terms $\alpha(\rho)$ and $\omega(\rho)$
of the extended normal form of the autonomous SSM, $\mathcal{M}_{0}$.
These terms are constructed from trajectories of the unforced system,
thus eq. \eqref{eq:closedformFRC} predicts the behavior of a non-linearizable
dynamical system under forcing based solely on unforced training data.
The stability of the predicted periodic response follows from a simple
linear analysis at the corresponding fixed point of the ODE \eqref{eq:periodicSSMforcing}.
The first formula in \eqref{eq:closedformFRC} also contains another
frequently used notion of nonlinear vibration analysis, the \emph{dissipative
backbone curve} $\omega(\rho)$, which describes the instantaneous
amplitude-frequency relation along freely decaying vibrations within
the SSM. 

As we will also show in examples, our unforced model-based predictions
for forced periodic response (see the Methods section \ref{subsec:Prediction-of-forced-response-methods})
are confirmed by numerical simulation or dedicated laboratory experiments
on forced systems.

\subsection{Examples}

\begin{figure}
\includegraphics[width=1\textwidth]{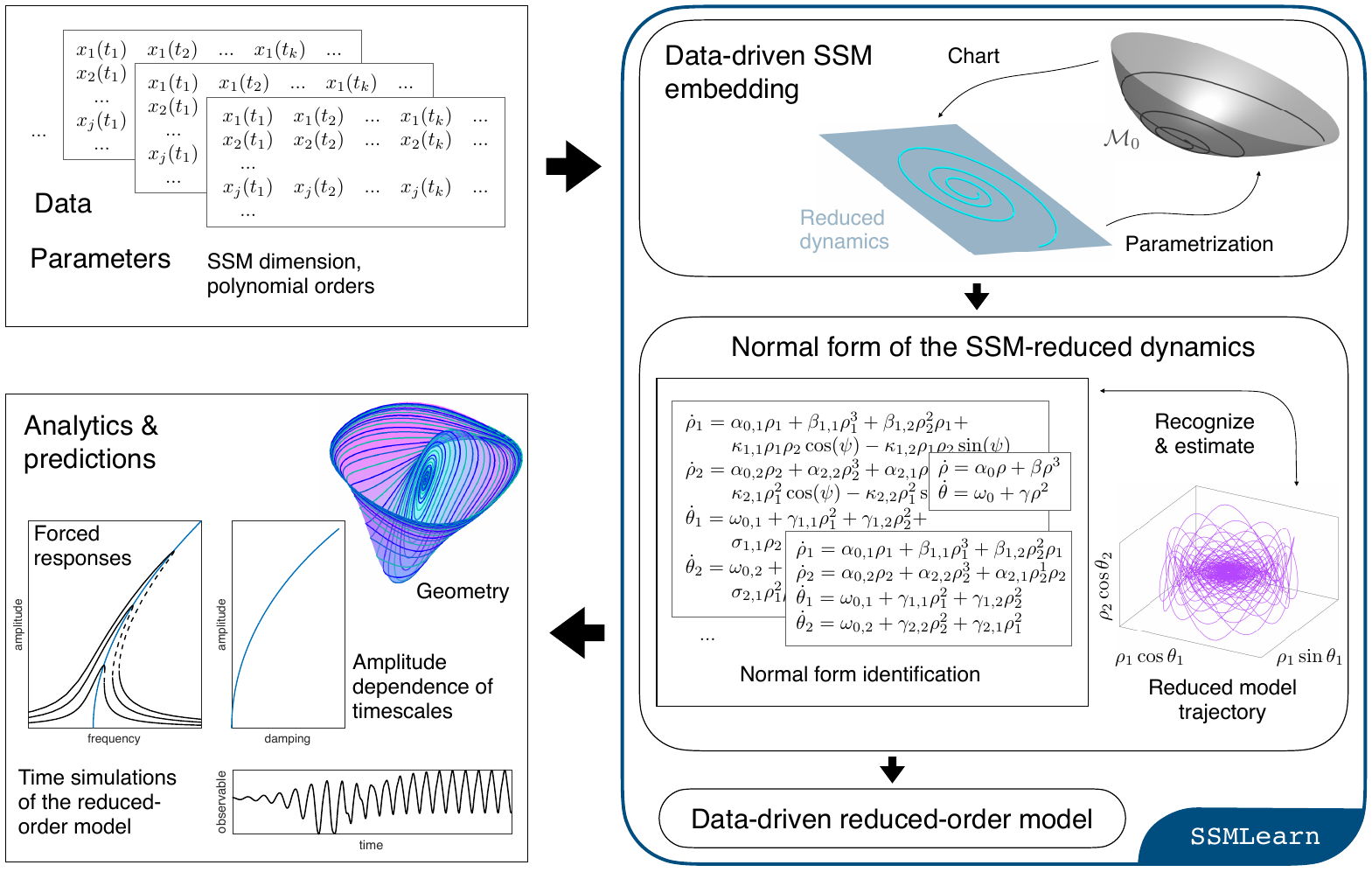}\caption{Schematics of the data-driven, SSM-based model reduction algorithm
implemented in \texttt{SSMLearn}.}\label{fig:ssmlearn}
\end{figure}

We now illustrate data-driven, SSM-based modeling and prediction on
several numerical and experimental data sets describing non-linearizable
physical systems. Further applications are described in \cite{cenedese21mech}.
Both the numerical and the experimental data sets were initialized
without knowledge of the exact SSM. All our computations have been
carried out by the publicly available \textsc{MATLAB}\textsuperscript{®}
package, \texttt{SSMLearn}, whose repository also contains further
examples not discussed here. The main algorithm behind \texttt{SSMLearn}
is illustrated in Fig. \ref{fig:ssmlearn}, with more detail given
in the Methods section \ref{subsec:Summary-of-the-algorithm}.

To quantify the errors of an SSM-based reduced model, we use the \emph{normalized
mean-trajectory-error} ($\mathrm{NMTE}$). For $P$ observations of
the observable vector $\mathbf{y}_{j}$ and their model-based reconstructions,
$\hat{\mathbf{y}}_{j}$, this modeling error is defined as 
\begin{equation}
\mathrm{NMTE}=\frac{1}{\|\underline{\mathbf{y}}\|}{\displaystyle \frac{1}{P}\sum_{j=1}^{P}\left\Vert \mathbf{y}_{j}-\hat{\mathbf{y}}_{j}\right\Vert }\,.\label{eq rRMSE}
\end{equation}
Here $\underline{\mathbf{y}}$ is a relevant normalization vector,
such as the data point with the largest norm. When validating the
reduced dynamics for a given testing trajectory, we run the reduced
model from the same initial condition for the comparison. Increasing
the order of the truncated normal form polynomials in eq. (\ref{eq:oscillatorsnormalform})
generally reduces the $\mathrm{\mathrm{NMTE}}$ error to any required
level but excessively small errors can lead to overfitting. In our
examples, we will be allowing model errors in the order of $1\%-4\%$
to avoid overfitting.

\subsubsection{Finite-element model of a damped-forced beam}

\begin{figure}[!t]
\includegraphics[width=1\textwidth]{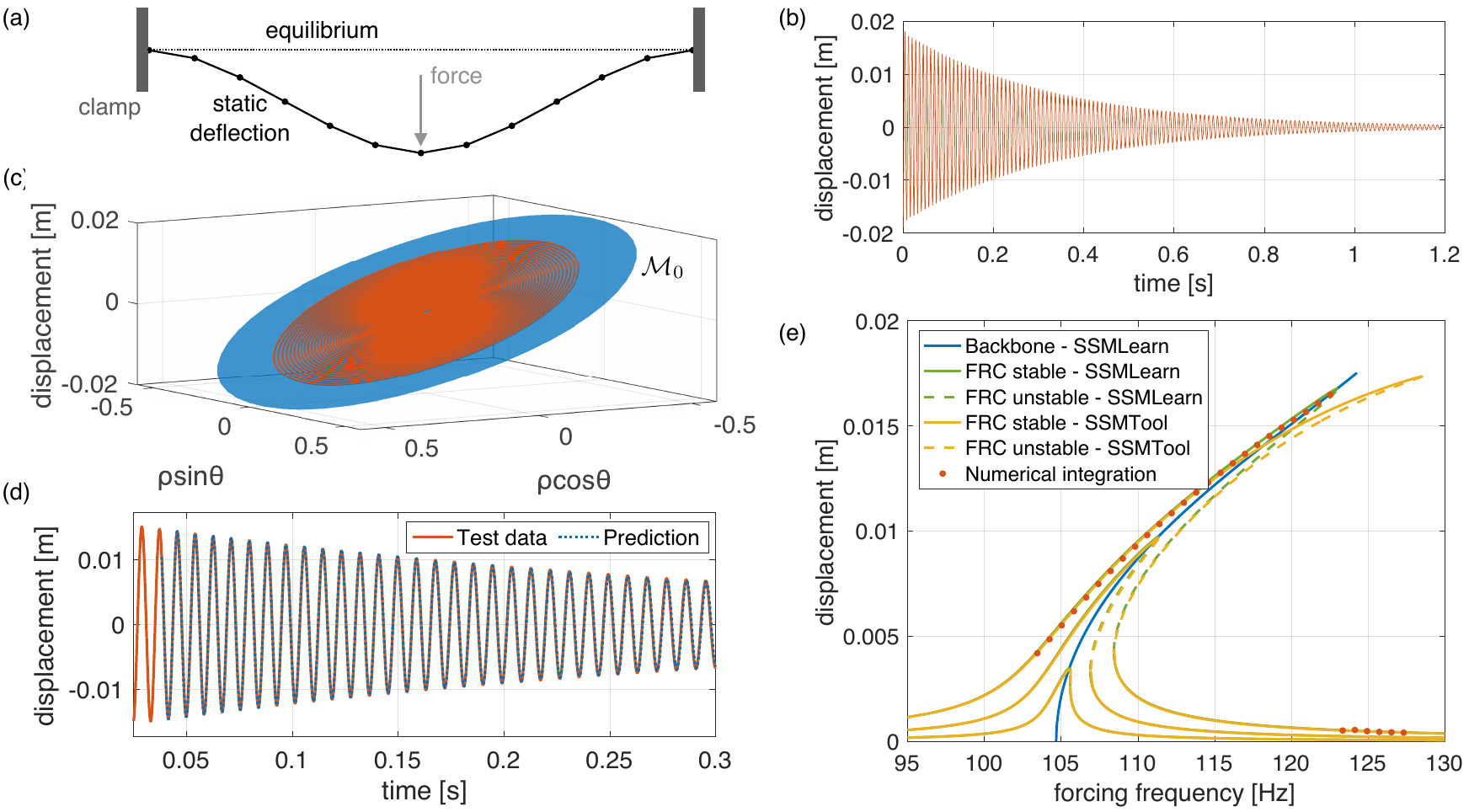}\caption{Construction of a data-driven nonlinear reduced-order model on the
slowest SSM of a von Kármán beam. (a) System setup and the initial
condition for the decaying training trajectory shown in (b) in terms
of the midpoint displacement. (c) The SSM, $\mathcal{M}_{0}$, in
the delay embedding space, shown along with the reconstructed test
trajectory in extended normal form coordinates. (d) Zoom of the prediction
of the reduced order model for the test trajectory not used in learning
$\mathcal{M}_{0}$. (e) Closed-form backbone curve and forced response
curve (FRC) predictions ($\epsilon>0$, $\ell=1$) by \texttt{SSMLearn}
are compared with analytic FRC calculations performed by \texttt{SSMTool}
\cite{jain2021} and with results from numerical integration of the
forced-damped beam.}\label{fig:vonkarmanbeam}
\end{figure}

We consider a finite-element discretization of a von Kármán beam with
clamped-clamped boundary conditions \cite{jain2018}, shown in Fig.
\ref{fig:vonkarmanbeam}(a). In contrast to the classic Euler-Bernoulli
beam, the von Kármán model captures moderate deformations by including
a nonlinear, quadratic term in the kinematics. We first construct
a 33 degree-of-freedom, damped, unforced finite element model (i.e.,
$n=66$ and $\epsilon=0$ in eq. \eqref{eq:1storder_system}) for
an aluminum beam of length 1 {[}m{]}, width 5 {[}cm{]}, thickness
2 {[}cm{]} and material damping modulus 10$^{6}$ {[} {]} (see the
Supplementary Information for more detail).

Our objective is to learn from numerically generated trajectory data
the reduced dynamics on the slowest, two-dimensional SSM, $W(E_{1})$,
of the system, defined over the slowest two-dimensional ($d=2$) eigenspace
$E_{1}$ of the linear part. To do so, we generate two trajectories
starting from initial beam deflections caused by static loading of
12 {[}kN{]} and 14 {[}kN{]} at the midpoint, as shown in Fig. \ref{fig:vonkarmanbeam}(a).
The latter trajectory, shown in Fig. \ref{fig:vonkarmanbeam}(b),
is used for training, the other for testing. Along the trajectories,
we select our single observable $s(t)$ to be the midpoint displacement
of the beam.

The beam equations are analytic $(r=a)$, and hence the SSM, $W(E_{1})$,
admits a convergent Taylor expansion near the origin. The minimal
embedding dimension for the two-dimensional, $W(E_{1})$, as required
by Whitney's theorem, is $p=5$, which is not satisfied by our single
scalar observable $s(t).$ We therefore employ delay-embedding using
$\mathbf{y}(t)=\left(s(t),s(t+\Delta t),\ldots,s(t+4\Delta t)\right)$
with $\Delta t=0.0955$ {[}ms{]}. By Takens's theorem, this delayed
observable embeds the SSM in $\mathbb{R}^{5}$ with probability one. 

A projection of the embedded SSM, $\mathcal{M}_{0}\in\mathbb{R}^{5},$
onto three coordinates is shown in Fig. \ref{fig:vonkarmanbeam}(c).
On $\mathcal{M}_{0}$, \texttt{SSMLearn} returns the $7^{th}$-order
extended normal form 
\begin{align}
\dot{\rho} & =\alpha(\rho)\rho,\quad\alpha(\rho)=-3.02-5.79\rho^{2}+57.5\rho^{4}-191\rho^{6},\nonumber \\
\dot{\theta} & =\omega(\rho),\,\,\,\quad\omega(\rho)=658+577\rho^{2}-347\rho^{4}-387\rho^{6},\label{eq:vonkarmanbeamROM}
\end{align}
to achieve our preset reconstruction error bar of $3\%$ on the test
trajectory ($\mathrm{NMTE}=0.027$), shown in Fig. \ref{fig:vonkarmanbeam}(d). 

We now use the model (\ref{eq:vonkarmanbeamROM}), trained on a single
decaying trajectory, to predict the forced response of the beam for
various forcing amplitudes and frequencies in closed form. We will
then compare these predictions with analytic forced response computations
for the forced SSM, $\mathcal{M}_{\epsilon}(\Omega t)$, obtained
from \texttt{SSMTool} \cite{jain2021} and with numerical simulations
of the damped-forced beam. The periodic forcing is applied at the
midpoint node; the Taylor expansion order in \texttt{SSMTool} for
the analytically computed dynamics on $\mathcal{M}_{\epsilon}(\Omega t)$
is set to $7$, as in (\ref{eq:vonkarmanbeamROM}). Figure \ref{fig:vonkarmanbeam}(e)
shows the FRCs (green) and the backbone curve (blue) predicted by
\texttt{SSMLearn} based on formula \eqref{eq:closedformFRC} from
the single unforced trajectory in Fig. \ref{fig:vonkarmanbeam}(b).
To obtain the relevant forcing amplitudes $f$ in the delay-observable
space, we have followed the calibration procedure described in the
Methods section \ref{subsec:Prediction-of-forced-response-methods}
for the forcing values $\left|\epsilon\mathbf{f}_{1}\right|=15,45,95$
{[}N{]} at the single forcing frequency $\Omega=103.5$ {[}Hz{]}.
Recall that coexisting stable (solid lines) and unstable (dashed lines)
periodic orbits along the same FRC are hallmarks of non-linearizable
dynamics and hence cannot be captured by the model reduction techniques
we reviewed in the Introduction for linearizable systems.

The data-based prediction for the FRCs agrees with the analytic FRCs
for low forcing amplitudes but departs from it for higher amplitudes.
Remarkably, as the numerical simulations (red) confirm, the data-based
FRC is the correct one. The discrepancy between the two FRCs for large
amplitudes only starts decreasing under substantially higher-order
Taylor series approximations used in \texttt{SSMTool} (see the Supplementary
Information). This suggests the use of the data-based approach for
this class of problems even if the exact equations of motion are available. 

\subsubsection{Vortex-shedding behind a cylinder\label{subsec:Vortex-shedding-behind-a-cylinder}}

\begin{figure}[!t]
\includegraphics[width=1\textwidth]{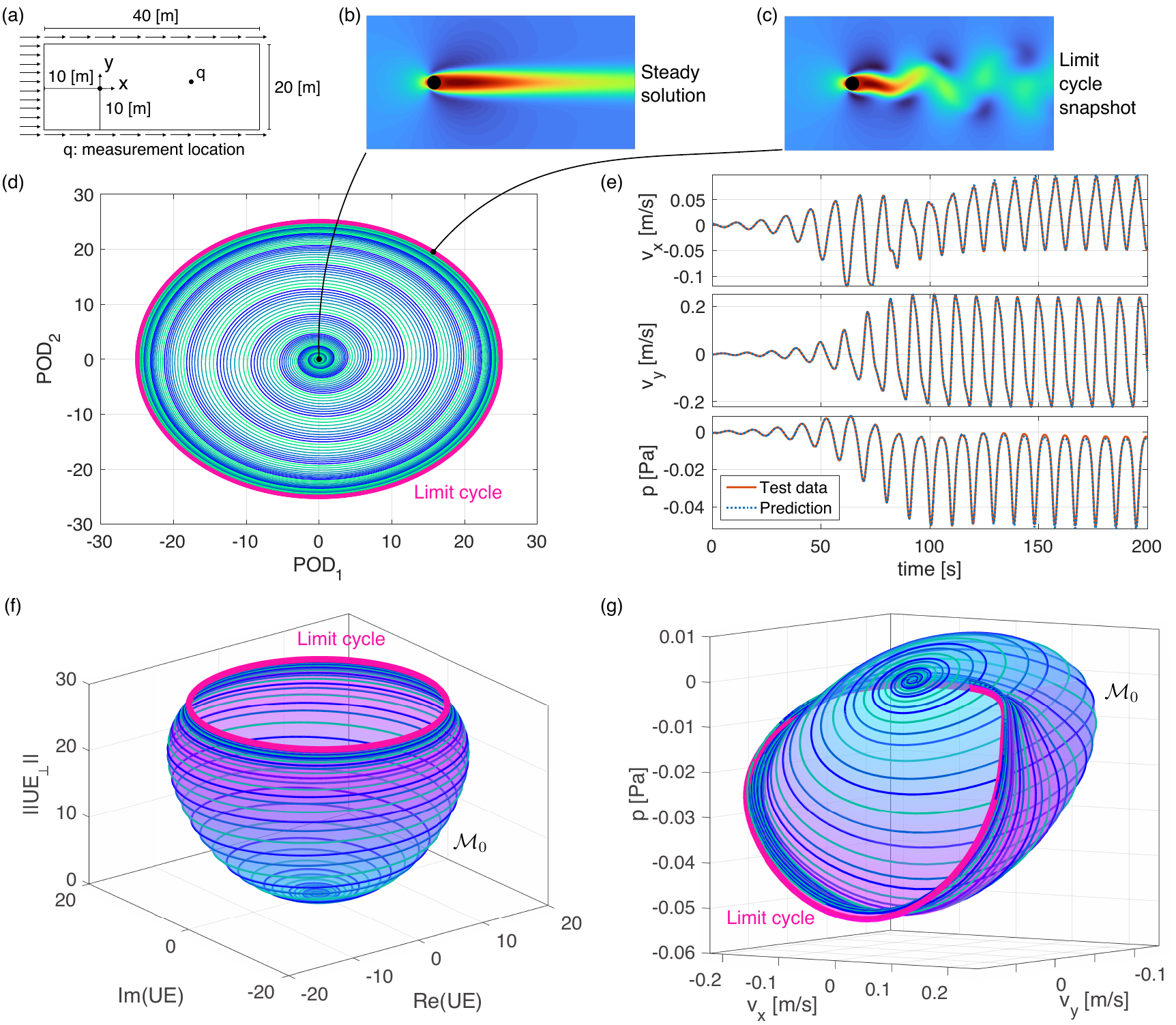}\caption{Data-driven nonlinear SSM-reduced model on the unstable manifold of
the steady solution of the flow past a cylinder. (a) Problem setup.
(b,c) Snapshots of the steady solution and the time-periodic vortex-shedding
solution (limit cycle, in magenta). (d) Trajectories projected on
the 2-dim. subspace spanned by the two-leading POD modes of the limit
cycle. (e) Model-based reconstruction of the test trajectory (not
used in learning the SSM) in terms of velocities and pressures measured
at a location $q$ shown in plot (a).  (f) The SSM formed by the unstable
manifold of the origin, along with some reduced trajectories, plotted
over the unstable eigenspace $UE\equiv E_{1}$; $\|\mathrm{UE}_{\perp}\|$
denotes the normed projection onto the orthogonal complement $\mathrm{UE}_{\perp}$.
(g) Same but projected over velocity and pressure coordinates.}\label{fig:vortexshedding}
\end{figure}

Next we consider the classic problem of vortex shedding behind a cylinder
\cite{loiseau20}. Our input data for SSM-based reduced modeling
are the velocity and pressure fields over a planar, open fluid domain
with a hole representing the cylinder section, as shown in Fig. \ref{fig:vortexshedding}(a).
The boundary conditions are no-slip on the circular inner boundary,
standard outflow on the outer boundary at the right side, and fixed
horizontal velocity on the three remaining sides \cite{barkley1996}.
The Reynolds number for this problem is the ratio between the cylinder
diameter times the inflow velocity and the kinematic viscosity of
the fluid. 

Available studies \cite{barkley1996,noack2003,loiseau20} report
that, at low Reynolds number, the two-dimensional unstable manifold,
$W^{u}(SS)$, of the wake-type steady solution, $SS$, in Fig. \ref{fig:vortexshedding}(b)
connects $SS$ to the limit cycle shown in Fig. \ref{fig:vortexshedding}(c).
Here we evaluate the performance of \texttt{SSMLearn} on learning
this unstable manifold as an SSM, along with its reduced dynamics,
from trajectory data at Reynolds number equal to 70. For this SSM,
we again have $d=2$ and $r=a$, as in our previous example. There
is no external forcing in this problem, and hence we have $\epsilon=0$
in eq. \eqref{eq:1storder_system}. In contrast to prior studies that
often consider a limited number of observables \cite{noack2003,rowley17,loiseau20},
here we select the full phase space of the discretized Navier-Stokes
simulation to be the observable space for illustration, which yields
$n=p=76,876$ in eq. \eqref{eq:1storder_system}. We generate nine
trajectories numerically, eight of which will be used for training
and one for testing the SSM-based model. 

The nine initial conditions of our input trajectory data are small
perturbations from the wake-type steady solution along its unstable
directions, equally spaced on a small amplitude circle on this unstable
plane. All nine trajectories quickly converge to the unstable manifold
and then to the limit cycle representing periodic vortex shedding.

We choose to parametrize the SSM, $\mathcal{\mathcal{M}}_{0}=W^{u}(SS)$,
with two leading POD modes of the limit cycle, which have been used
in earlier studies for this problem. The training trajectories projected
onto these two POD modes are shown in Fig. \ref{fig:vortexshedding}(d).
To limit the modeling error \eqref{eq rRMSE} to less than $\mathrm{NMTE}=1$
$\%$, \texttt{SSMLearn} requires a polynomial order of 18 in the
SSM computations. For this order, our approach can accommodate the
strong mode deformation observed for this problem \cite{noack2003},
manifested by a fold of the SSM over the unstable eigenspace in Fig.
\ref{fig:vortexshedding}(f). Figure \ref{fig:vortexshedding}(g)
shows the strongly nonlinear geometry of $\mathcal{\mathcal{M}}_{0}$
projected to the observable subspace formed by the velocities and
the pressure of a probe point in the wake.

To capture the SSM-reduced dynamics with acceptable accuracy, we need
to compute the extended normal form up to order 11, obtaining
\begin{equation}
\begin{array}{rl}
\dot{\rho}=\alpha(\rho)\rho= & 0.0584\rho-0.479\rho^{3}+1.27\rho^{5}+6.80\rho^{7}-58.9\rho^{9}+108\rho^{11},\\
\dot{\theta}=\omega(\rho)\,\,\,= & 0.553+0.441\rho^{2}-3.38\rho^{4}+55.5\rho^{6}-321\rho^{8}+626\rho^{10}.
\end{array}\label{eq:vortexsheddingROM}
\end{equation}
To describe a transition qualitatively from a fixed point to a limit
cycle, the reduced-order dynamical model should be at least of cubic
order \cite{noack2003}. Capturing the qualitative behavior (i.e.,
the unstable fixed point and the stable limit cycle), however, does
not imply a low $\mathrm{NMTE}$ error for the model. Indeed, the
data-driven cubic normal form for this example gives a reconstruction
error of $\mathrm{NMTE}=117$ $\%$ normalized over the limit cycle
amplitude, mainly arising from an out-of-phase convergence to the
limit cycle along the testing trajectory. In contrast, the $\mathcal{O}(11)$
normal form in eq. (\ref{eq:vortexsheddingROM}) reduced this error
drastically to $\mathrm{NMTE}=3.86$ $\%$ on the testing trajectory,
as shown in Fig. \ref{fig:vortexshedding}(e).

We show in Section 1.2.3 of the Supplementary Information that
for comparable accuracy, the Sparse Identification of Nonlinear DYnamics
(SINDy) approach of \cite{brunton2016} returns non-sparse nonlinear
models for this example. Similarly, while the DMD \cite{kutz16}
can achieve highly accurate curve-fitting on the available training
trajectories with a high-dimensional linear model, that model only
captures linearizable dynamics near the origin. As a consequence,
its trajectories grow without bound over longer integration times
and hence fail to capture the limit cycle.

\subsubsection{Fluid sloshing experiments in a tank\label{subsec:Fluid-sloshing-experiments}}

Fluid oscillations in a tank exhibit highly nonlinear characteristics
\cite{taylor1953}. To describe such non-linearizable softening effects
observed in the sloshing motion of surface waves, Duffing-type models
have been proposed \cite{ockendon1973}. While amplitude variations
observed in forced experiments can be fitted to forced softening Duffing
equations, nonlinear damping remains a challenge to identify \cite{bauerlein2021}.
\begin{figure}
\includegraphics[width=1\textwidth]{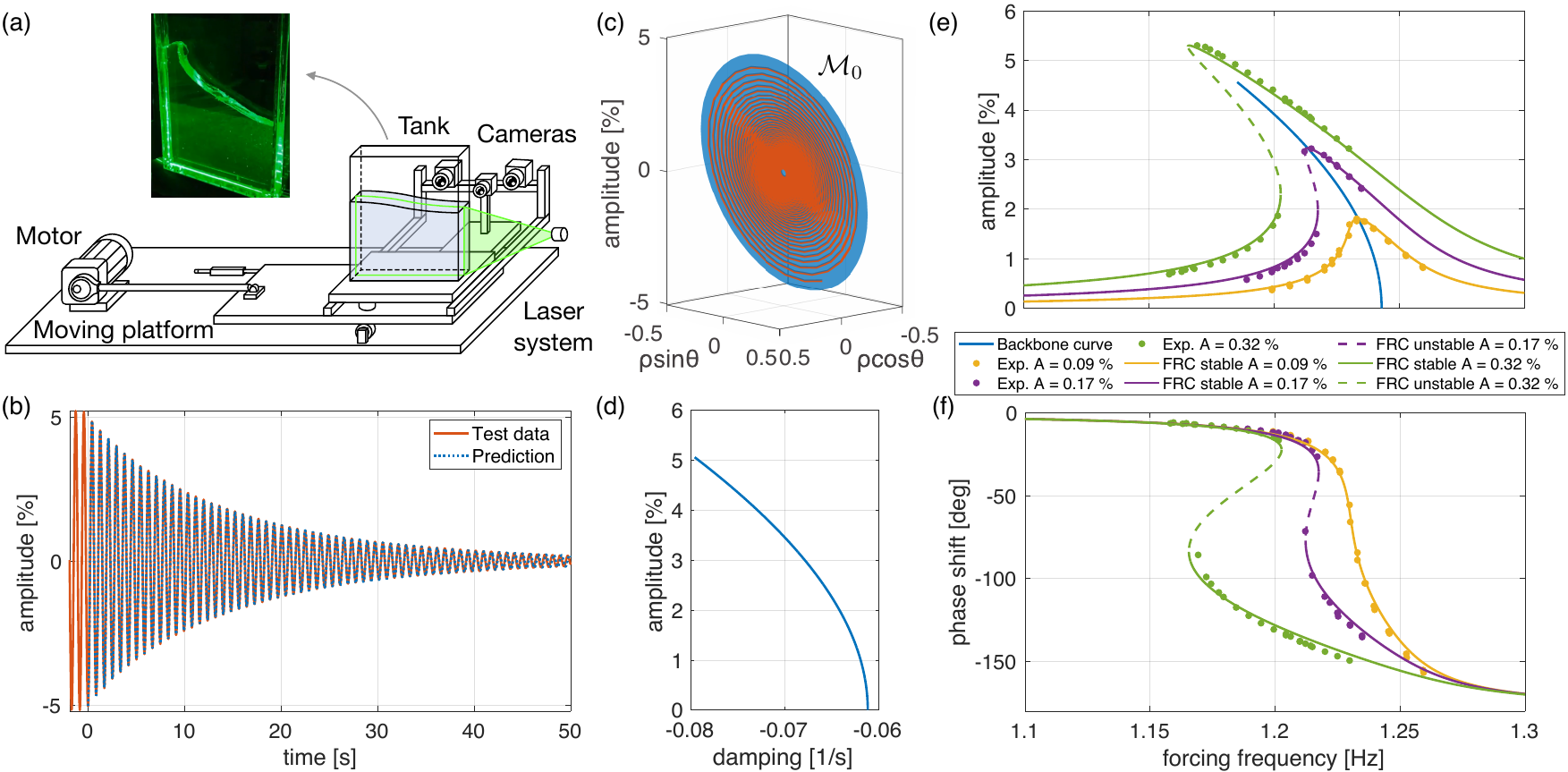}\caption{Data-driven nonlinear reduced-order model on the slowest SSM of fluid
sloshing in a tank. (a) Setup for the sloshing experiment \cite{bauerlein2021}.
(b) Decaying model-testing trajectory and its reconstruction from
an unforced, SSM-based model (c) The geometry of the embedded SSM
(d) Nonlinear damping $\alpha(\rho)$ from the SSM-reduced dynamics
(e,f) Closed form, SSM-based predictions of the FRCs and the response
phases $\psi_{0}$ for three different forcing amplitudes (solid lines),
with their experimental confirmation superimposed (dots).}\label{fig:sloshing}
\end{figure}

The experiments we use to construct an SSM-reduced nonlinear model
for sloshing were performed in a rectangular tank of width 500 {[}mm{]}
and depth 50 {[}mm{]}, partially filled with water up to a height
of 400 {[}mm{]}, as shown in Fig. \ref{fig:sloshing}(a). The tank
was mounted on a platform excited harmonically by a motor. The surface
level was detected via image processing from a monochrome camera.
As an observable\textbf{ }$s(t)$ we used the horizontal position
of the computed center of mass of the water at each time instant,
normalized by the tank width. This physically meaningful scalar is
robust with respect to image evaluation errors \cite{bauerlein2021}.

We identify the unforced nonlinear behavior of the system from data
obtained in resonance decay experiments \cite{peeters2011a}. In
those experiments (as in Fig. \ref{fig:sloshing}a, but with a shaker
instead of a motor), once a periodic steady state is reached under
periodic horizontal shaking of the tank, the shaker is turned off
and the decaying sloshing is recorded. We show such a decaying observable
trajectory (orange line) in Fig. \ref{fig:sloshing}(b), with the
shaker switched off slightly before zero time. This damped oscillation
is close, by construction, to the two-dimensional, slowest SSM of
the system. We use three such decaying observer trajectories (two
for training and one for model testing) for the construction of a
two-dimensional ($d=2$), autonomous, SSM-based reduced order model
for $s(t)$. For delay embedding dimension, we again pick $p=5$,
the minimal value guaranteed to be generically correct for embedding
the SSM by Takens's theorem. The delay used in sampling $s(t)$ is
$\Delta t=0.033$ {[}s{]}. For this input and for a maximal reconstruction
error of $2\%$, \texttt{SSMLearn} identifies a nearly flat SSM in
the delayed observable space--see Fig. \ref{fig:sloshing}(c)--with
a cubic extended normal form 
\begin{equation}
\dot{\rho}=-0.0628\rho-0.0572\rho^{3},\quad\dot{\theta}=7.80-1.67\rho^{2}.\label{eq:sloshingROM}
\end{equation}
This lowest-order, Stuart--Landau-type normal form (see \eqref{eq:cubicnormalform})
already constitutes an accurate reduced order model with $\mathrm{NMTE}=1.88\%$
on the testing data set (see Fig. \ref{fig:sloshing}(b)). The amplitude-dependent
nonlinear damping, $\alpha(\rho)$, provided by this model is plotted
in Fig. \ref{fig:sloshing}(d) with respect to the physical amplitude. 

In another set of experiments with the setup of Fig. \ref{fig:sloshing}a,
steady states of periodically forced sloshing were measured in sweeps
over a range of forcing frequencies under three different shaker amplitudes.
As in the previous beam example, we identify the corresponding forcing
amplitude, $f$, in \eqref{eq:periodicSSMforcing} at the maximal
amplitude response of each frequency sweep. Shown in Fig. \ref{fig:sloshing}(e,f),
the closed-form predictions for FRCs from eq. \eqref{eq:closedformFRC}
(solid lines) match closely the experimental FRCs (dots). Given the
strong nonlinearity of the FRC, any prediction of this curve from
a DMD-based model is bound to be vastly inaccurate, as we indeed show
in Section 1.3 of the Supplementary information. 

The phase $\psi_{0}$ of the forced response relative to the forcing
has been found difficult to fit to forced Duffing-type models \cite{bauerlein2021},
but the present modeling methodology also predicts this phase accurately
using the second expression in \eqref{eq:closedformFRC}. The blue
curve in Fig. \ref{fig:sloshing}(d) shows the backbone curve of decaying
vibrations, which terminates at the highest amplitude occurring in
the training data set. This plot therefore shows that the closed-form
FRC predictions obtained from the SSM-based reduced model are also
effective for response amplitudes outside the training range of the
reduced model.

\section{Discussion}

We have described a data-driven model reduction procedure for non-linearizable
dynamical systems with coexisting isolated stationary states. Our
approach is based on the recent theory of spectral submanifolds (SSMs),
which are the smoothest nonlinear continuations of spectral subspaces
of the linearized dynamics. Slow SSMs form a nested hierarchy of attractors
and hence the dynamics on them provide a hierarchy of reduced order
models with which generic trajectories synchronize exponentially fast.
These SSMs and their reduced models smoothly persist under moderate
external forcing, yielding low-dimensional, mathematically exact reduced
order models for forced versions of the same dynamical system. The
normal hyperbolicity of SSMs also ensures their robustness under small
noise.

All these results have been implemented in the open source \textsc{MATLAB}\textsuperscript{®}
package, \texttt{SSMLearn}, which we have illustrated on data sets
arising from forced nonlinear beam oscillations, vortex shedding behind
a cylinder and water sloshing in a vibrating tank. For all three examples,
we have found that two-dimensional data-driven extended normal forms
on the slowest SSMs provide sparse yet accurate models of non-linearizable
dynamics in the space of the chosen observables. Beyond matching training
and testing data, SSM-reduced models prove their intrinsic, qualitative
meaning by predicting non-linearizable, forced steady states purely
from decaying, unforced data. 

In this brief report, examples of higher-dimensional SSMs and multi-harmonic
forcing have not been considered, even though \texttt{SSMLearn} is
equipped to handle them. Higher-dimensional SSMs are required in the
presence of internal resonances or in non-resonant problems in which
initial transients also need to be captured more accurately. A limitation
of our approach for non-autonomous systems is the assumption of quasiperiodic
external forcing. Note, however, that even specific realizations of
stochastic forcing signals can be approximated arbitrarily closely
with quasiperiodic functions over any finite time interval of interest.
A further limitation in our work is the assumption of smooth system
dynamics. For data from non-smooth systems, \texttt{SSMLearn} will
nevertheless return an equivalent smooth reduced-order model whose
accuracy is a priori known from the available mean-squared error of
the SSM fitting and conjugacy error of the normal form construction.
We are addressing these challenges in ongoing work to be reported
elsewhere. Further applications of \texttt{SSMLearn} to physical problems
including higher-dimensional coexisting steady states (see, e.g.,
\cite{deng20}) are also underway.

\appendix

\section{Methods}

\subsubsection*{Data availability}

All data discussed in the results presented here is publicly
available in the \texttt{SSMLearn} repository at \url{github.com/haller-group/SSMLearn}.

\subsubsection*{Code availability}

The code supporting the results presented here is publicly
available in the \texttt{SSMLearn} repository at \url{github.com/haller-group/SSMLearn}.

\subsection{Existence of SSMs\label{subsec:Existence-of-SSMs}}

In the context of rigid body dynamics, invariant manifolds providing
generalizations of invariant spectral subspaces to nonlinear systems
were first envisioned and formally constructed as nonlinear normal
modes by \cite{shaw93} (see \cite{renson16} for a recent review
of related work). Later studies, however, pointed out the non-uniqueness
of nonlinear normal modes in specific examples (\cite{neild15,cirillo16}). 

In the mathematics literature, \cite{cabre03} obtained general results
on the existence, smoothness and degree of uniqueness of such invariant
manifolds for mappings on Banach spaces. These results use a special
parameterization method to construct the manifolds even in evolutionary
partial differential equations that admit a well-posed flow map in
both time directions (see \cite{kogelbauer18} for a mechanics application).
The results have been extended to a form applicable to dynamical systems
with quasiperiodic time dependence \cite{haro06}. An extensive account
of the numerical implementation of the parametrization method with
a focus on computing invariant tori and their whiskers in Hamiltonian
systems is also available \cite{haro16}.

\cite{haller16} discussed the existence of the SSM, $W(E,\boldsymbol{\Omega}t;\epsilon)$,
depending on its \emph{absolute spectral quotient},
\begin{equation}
\Sigma(E)=\mathrm{Int}\,\left[\frac{{\displaystyle \max_{\lambda\in\mathrm{Spect}(\mathbf{A}\vert_{S})}|\mathrm{Re}\lambda}|}{{\displaystyle \min_{\lambda_{e}\in\mathrm{Spect}(\mathbf{A}\vert_{E})}|\mathrm{Re}\lambda_{e}|}}\right],\label{eq:absolute spectral quotient}
\end{equation}
where $\mathrm{Spect}(\mathbf{A}\vert_{S})$ is the stable (unstable)
spectrum of $\mathbf{A}$ if the SSM is stable (unstable). For a stable
SSM, $\Sigma(E)$ is the integer part of the quotient of the minimal
real part in the spectrum of $\mathbf{A}$ and the maximal real part
of the spectrum of $\mathbf{A}$ restricted to $E$. 

Based on $\Sigma(E)$, we call a $d$-dimensional spectral subspace
$E$ \emph{non-resonant} if for any set $\left(m_{1},\ldots,m_{d}\right)$
of nonnegative integers satisfying $2\leq\sum_{j=1}^{d}m_{j}\leq\Sigma(E)$,
the eigenvalues, $\lambda_{k}$, of $\mathbf{A}$ satisfy
\begin{equation}
\sum_{j=1}^{d}m_{j}\mathrm{Re}\lambda_{j}\neq\mathrm{Re}\lambda_{k},\quad\lambda_{k}\in\mathrm{Spect}(\mathbf{A})-\mathrm{Spect}(\mathbf{A}\vert_{E}).\label{eq:forced_nonresonance}
\end{equation}
This condition only needs to be verified for resonance orders between
$2$ and $\Sigma(E)$ \cite{haro06}. In particular, a $1\colon1$
resonance between $E_{1}$ and $E_{2}$ is allowed if $\dim E_{1}=\dim E_{2}=1$,
in which case each strongly resonant spectral subspace gives rise
to a unique nearby spectral submanifold. 

If $E$ violates the non-resonance condition \eqref{eq:forced_nonresonance},
then $E$ can be enlarged to a higher-dimensional spectral subspace
until the non-resonance relationship \eqref{eq:forced_nonresonance}
is satisfied. In the absence of external forcing ($\epsilon=0$),
the non-resonance condition \eqref{eq:forced_nonresonance} can also
be relaxed with the help of the \emph{relative spectral quotient},
\begin{equation}
\sigma(E)=\mathrm{Int}\,\left[\frac{{\displaystyle \max_{\lambda\in\mathrm{Spect}(\mathbf{A}\vert_{S})-\mathrm{Spect}(\mathbf{A}\vert_{E})}|\mathrm{Re}\lambda}|}{{\displaystyle \min_{\lambda_{e}\in\mathrm{Spect}(\mathbf{A}\vert_{E})}|\mathrm{Re}\lambda_{e}|}}\right],\label{eq:relative spectral quotient}
\end{equation}
to the form
\begin{equation}
\sum_{j=1}^{d}m_{j}\lambda_{j}\neq\lambda_{k},\quad\lambda_{k}\in\mathrm{Spect}(\mathbf{A})-\mathrm{Spect}(\mathbf{A}\vert_{E}),\qquad2\leq\sum_{j=1}^{d}m_{j}\leq\sigma(E).\label{eq:unforced_nonresonance}
\end{equation}
This is indeed a relaxation because condition \eqref{eq:unforced_nonresonance}
is only violated if both the real and the imaginary parts of eigenvalues
involved are in the exact same resonance with each other. In contrast,
\eqref{eq:forced_nonresonance} is already violated when the real
parts are in resonance with each other. 

If $\mathrm{Re}\lambda_{1}<0$ in eq. \eqref{eq:eigenvalue_ordering}
and all $E^{k}$ subspaces are non-resonant, then the nested set of
slow spectral submanifolds,
\[
W(E^{1},\boldsymbol{\Omega}t;\epsilon)\subset W(E^{2},\boldsymbol{\Omega}t;\epsilon)\subset W(E^{3},\boldsymbol{\Omega}t;\epsilon)\subset\ldots,
\]
gives a hierarchy of local attractors. All solutions in a vicinity
of $\mathbf{x}=\mathbf{0}$ approach the reduced dynamics on one of
these attractors exponentially fast, as sketched in Fig. \ref{fig:linear vs nonlinear dynamics}b
for the $\epsilon=0$ limit. As we will see, non-linearizable dynamics
tend to emerge on $W(E^{k},\boldsymbol{\Omega}t;\epsilon)$ due to
near-resonance between the linearized frequencies within $E^{k}$
and the forcing frequencies $\boldsymbol{\Omega}$. The specific location
of nontrivial steady states in $W(E^{k},\boldsymbol{\Omega}t;\epsilon)$
is then determined by a balance between the nonlinearities, damping
and forcing. 

A resonant $E^{k}$ subspace can be enlarged by adding the next $k'$
modal subspaces to it until $E^{k+k'}$ in the hierarchy \eqref{eq:hierarchy of slow spectral subspaces}
becomes non-resonant and hence admits an SSM, $W(E^{k+k'},\boldsymbol{\Omega}t;\epsilon)$.
This technical enlargement is also in agreement with the physical
expectation that all interacting modes have to be included in an accurate
reduced-order model. Finally, we note that SSMs are robust features
of dynamical systems: they inherit smooth dependence of the vector
field in (\ref{eq:1storder_system}) on parameters \cite{haller16}. 

For discrete-time dynamical systems of the form
\begin{equation}
\mathbf{x}_{k+1}=\tilde{\mathbf{A}}\mathbf{x}_{k}+\tilde{\mathbf{f}}_{0}(\mathbf{x}_{k})+\epsilon\tilde{\mathbf{f}}_{1}(\mathbf{x}_{k},\boldsymbol{\phi}_{k};\epsilon),\qquad\boldsymbol{\phi}_{k+1}=\boldsymbol{\phi}_{k}+\tilde{\boldsymbol{\Omega}},\label{eq:1storder_map}
\end{equation}
the above results on SSMs apply based on the eigenvalues $\mu_{k}$
of $\tilde{\mathbf{A}}$. One simply needs to replace $\lambda_{k}$
with $\log\mu{}_{k}$ and $\mathrm{Re}\lambda_{k}$ with $\log|\mu{}_{k}|$
in formulas \eqref{eq:absolute spectral quotient}-\eqref{eq:unforced_nonresonance}
\cite{haller16}.

We close by noting that in a neighborhood of an SSM, an invariant
family of surfaces resembling the role of coordinate planes in a linear
system exists \cite{szalai20}. This invariant spectral foliation
(ISF) can, in principle, be used to generate a nonlinear analogue
of linear modal superposition in a vicinity of a fixed point. Constructing
the ISF from data has shown both initial promise and challenges to
be addressed.

\subsection{Embedding the SSM in the observable space\label{subsec:SSM-geometry-in-observable-space}}

Originally conceived for autonomous systems, the Takens delay embedding
theorem \cite{takens81} has been strengthened and generalized to
externally forced dynamics \cite{stark99}. By these results, the
embedding for a $d$-dimensional compact SSM subset, $\mathcal{C}\subset W(E,\boldsymbol{\Omega}t;\epsilon)$,
in the delay observable space as $\mathcal{M}(\boldsymbol{\Omega}t)$
is guaranteed for almost all choices of the observable $s(t)$ if
$p>2(d+l)$ and some generic assumptions regarding periodic motions
on $\mathcal{M}(\boldsymbol{\Omega}t)$ are satisfied \cite{sauer1991}. 

Of highest importance in technological applications is the case of
time-periodic forcing ($\ell=1$), with frequency $\boldsymbol{\Omega}=\Omega\in\mathbb{R}$
and period $T=2\pi/\Omega$. In this case, the Whitney and Takens
embedding theorems can be applied to the associated period-$T$ sampling
map (or Poincaré map) $\mathbf{P}_{t_{0}}\colon\mathbb{R}^{n}\to\mathbb{R}^{n}$
of the system based at time $t_{0}$. This map is autonomous and has
a time-independent SSM that coincides with the $d$-dimensional SSM,
$\mathcal{M}(\Omega t_{0})$, of the full system (\ref{eq:1storder_system}).
In this case, by direct application of the embedding theorems to the
discrete dynamical system generated by $\mathbf{P}_{t_{0}}$, the
typically sufficient embedding dimension estimate is improved to $p>2d$
for Whitney's and Takens's theorem.

Technically speaking, the available data will never be exactly on
an SSM, as these embedding theorems assume. By the smoothness of the
embeddings, however, points close enough to the SSM in the phase space
will be close to $\mathcal{M}(\boldsymbol{\Omega}t)$ in the observable
space under the embeddings. Moreover, as slow SSMs attract nearby
trajectories exponentially, the distance of observable data from the
embedded slow SSM will shrink exponentially fast. Therefore, even
under uncorrelated noise in the measurements, mean-squared estimators
are suitable for learning slow SSMs from data in the observable space,
as we illustrate in the Supplementary Information.

After a possible coordinate shift, the trivial fixed point of the
autonomous limit of system \eqref{eq:1storder_system} will be mapped
into the $\mathbf{y}=\mathbf{0}$ origin of the observable space.
To find an embedded, $d$-dimensional SSM, $\mathcal{\mathcal{M}}_{0}\in\mathbb{R}^{p}$
, attached to this origin for $\epsilon=0$, we focus on observable
domains in which $\mathcal{\mathcal{M}}_{0}$ is a graph over its
tangent space $T_{\mathbf{0}}\mathcal{M}_{0}$ at the origin $\mathbf{y}=\mathbf{0}$.
Such domains always exist and are generally large enough to capture
non-linearizable dynamics in most applications (but see below). Note
that $T_{\mathbf{0}}\mathcal{M}_{0}$ coincides with the image of
the spectral subspace $E$ in the observable space.

To learn such a graph-style parametrization for $\mathcal{\mathcal{M}}_{0}$
from data, we define a matrix $\mathbf{U}_{1}\in\mathbb{R}^{n\times d}$
with columns that are orthonormal vectors spanning the yet unknown
$T_{\mathbf{0}}\mathcal{M}_{0}$. The reduced coordinates $\boldsymbol{\eta}\in\mathbb{R}^{d}$
for a point $\mathbf{y}\in\mathcal{M}_{0}$ are then defined as the
orthogonal projection $\boldsymbol{\eta}=\mathbf{U}_{1}^{T}\mathbf{y}$.
We week a Taylor-expansion for $\mathcal{\mathcal{M}}_{0}$ near the
$\boldsymbol{\eta}=\mathbf{0}$ origin, denoting by $\boldsymbol{\eta}^{2:M}$
the family of all monomials of $d$ variables from degree 2 to $M$.
For example, if $d=2$ and $M=3$, then $\boldsymbol{\eta}^{2:3}=(\eta_{1}^{2},\eta_{1}\eta_{2},\eta_{2}^{2},\eta_{1}^{3},\eta_{1}^{2}\eta_{2},\eta_{1}\eta_{2}^{2},\eta_{2}^{3})^{T}$.
As a graph over $T_{\mathbf{0}}\mathcal{M}_{0}$, the manifold $\mathcal{\mathcal{M}}_{0}$
is approximated as $\mathbf{y}=\mathbf{V}_{1}\boldsymbol{\eta}+\mathbf{V}\boldsymbol{\eta}^{2:M}$,
where the matrices $\mathbf{V}_{1}$ and $\mathbf{V}$ contain coefficients
for the $d$-variate linear and nonlinear monomials, respectively.
Learning $\mathcal{\mathcal{M}}_{0}$ from a data set of $P$ observations
$\mathbf{y}_{1},\ldots,\mathbf{y}_{P}$ then amounts to finding the
$(\mathbf{U}_{1}^{*},\mathbf{V}_{1}^{*},\mathbf{V}^{*})$ matrices
that minimize the mean-square reconstruction error along the training
data:
\begin{equation}
\begin{split}(\mathbf{U}_{1}^{*},\mathbf{V}_{1}^{*},\mathbf{V}^{*})= & {\displaystyle \mathrm{arg}\min_{\mathbf{U}_{1},\mathbf{V}_{1},\mathbf{V}}\sum_{j=1}^{P}\left\Vert \mathbf{y}_{j}-\mathbf{V}_{1}\mathbf{U}_{1}^{T}\mathbf{y}_{j}-\mathbf{V}(\mathbf{U}_{1}^{T}\mathbf{y}_{j})^{2:M}\right\Vert ^{2},}\\
 & \mathbf{U}_{1}^{T}\mathbf{U}_{1}=\mathbf{I}.
\end{split}
\label{eq:learningparametrization}
\end{equation}
The simplest solution to this problem is $\mathbf{U}_{1}=\mathbf{V}_{1}$
with the additional constraint $\mathbf{V}_{1}^{T}\mathbf{V}=\mathbf{0}$,
which represents a basic nonlinear extension of the principal component
analysis \cite{bishop2006}. 

The above graph-style parametrization of the SSM breaks down for larger
$\mathbf{y}$ values if $\mathcal{\mathcal{M}}_{0}$ develops a fold
over $T_{0}\mathcal{M}_{0}$. That creates an issue for model reduction
if a nontrivial steady state on $\mathcal{\mathcal{M}}_{0}$ falls
outside the fold, as the limit cycle does in our vortex shedding example.
In that case, alternative parametrization methods for $\mathcal{\mathcal{M}}_{0}$
can be used to enhance the domain of the SSM-reduced model. These
methods include selecting the columns of $\mathbf{U}_{1}$ to be the
leading POD modes of the nontrivial steady state, or enlarging the
embedding space with (further) delayed observations. In these cases,
the columns of $\mathbf{V}_{1}$ are still orthonormal vectors spanning
$T_{\mathbf{0}}\mathcal{M}_{0}.$ 

In both Fig. \ref{fig:vonkarmanbeam} and Fig. \ref{fig:vonkarmanbeam},
the SSM, $\mathcal{M}_{0}$, is nearly flat in the delay-embedding
space. This turns out to be a universal property of delay embedding
for small delays and low embedding dimensions (see the Supplementary
Information).

For $\epsilon>0$ small (i.e., for moderate forcing), the autonomous
SSM, $\mathcal{M}_{0}$, already captures the bulk nonlinear behavior
of system \eqref{eq:1storder_system}. Indeed, for this forcing range,
the reduced dynamics on the corresponding SSM can simply be computed
as an additive perturbation of the autonomous dynamics on $\mathcal{M}_{0}$
\cite{breunung2018,ponsioen2019,ponsioen2020} (see section \ref{subsec:Prediction-of-forced-response}).

\subsection{SSM dynamics via extended normal forms\label{subsec:SSM-dynamics-via-extended normal forms}}

For an autonomous SSM $\mathcal{M}_{0}$, the reduced dynamics is
governed by a vector field 
\begin{equation}
\dot{\boldsymbol{\eta}}=\mathbf{r}(\boldsymbol{\eta})\label{eq:reduced dynamics}
\end{equation}
with a flow map $\boldsymbol{\varphi}_{\mathbf{r}}^{t}(\boldsymbol{\eta})$.
We can generically assume that the Jacobian $D\mathbf{r}(\boldsymbol{0})$
is semisimple, i.e., $D\mathbf{r}(\boldsymbol{0})\mathbf{B}=\mathbf{B}\boldsymbol{\Lambda}$,
where $\boldsymbol{\Lambda}\in\mathbb{C}^{d\times d}$ is a diagonal
matrix containing the eigenvalues of $D\mathbf{r}(\boldsymbol{0})$.
Classic normal form theory would seek to simplify the reduced dynamics
\eqref{eq:reduced dynamics} in a vicinity of $\boldsymbol{\eta}=\mathbf{0}$
via a nonlinear change of coordinates, $\boldsymbol{\eta}=\mathbf{h}(\mathbf{z})$,
so that the transformed vector field $\dot{\mathbf{z}}=\mathbf{n}(\mathbf{z})$
with flow map $\boldsymbol{\varphi}_{\mathbf{n}}^{t}(\mathbf{z})$
has a diagonal linear part and has as few nonlinear terms in its Taylor
expansion as possible. In our present setting, the origin is assumed
hyperbolic, in which case the classic normal form is simply $\dot{\mathbf{z}}=D\mathbf{r}(\boldsymbol{0})\mathbf{z}$
under appropriate non-resonance conditions that generically hold \cite{sternberg1958}.
The corresponding normal form transformation $\mathbf{h}(\mathbf{z})$,
however, is only valid on a small enough domain in which the dynamics
is linearizable.

To capture non-linearizable behavior, we employ extended normal forms
motivated by those used to unfold bifurcations \cite{guckenheimer83}.
In this approach, we construct normal forms that do not remove those
polynomial terms from \eqref{eq:reduced dynamics} whose removal would
result in small denominators in the Taylor coefficients $\mathbf{h}(\mathbf{z})$
and hence decrease its domain of convergence. Instead, we seek a normal
form for \eqref{eq:reduced dynamics} of the form
\begin{equation}
\begin{array}{c}
{\displaystyle \mathbf{n}(\mathbf{z};\mathbf{N})=\boldsymbol{\Lambda}\mathbf{z}+\mathbf{N}\mathbf{z}^{2:N},}\\
\mathbf{h}(\mathbf{z};\mathbf{H})=\mathbf{B}(\mathbf{z}+\mathbf{H}\mathbf{z}^{2:N}),\,\,\,\,\,\,\mathbf{h}^{-1}(\boldsymbol{\eta};\mathbf{H_{\star}})=\mathbf{B}^{-1}\boldsymbol{\eta}+\mathbf{H}_{\star}(\mathbf{B}^{-1}\boldsymbol{\eta})^{2:N},
\end{array}\label{eq:normalformSSMs}
\end{equation}
where the matrices $\mathbf{N}$, $\mathbf{H}$ and $\mathbf{H}_{\star}$
contain the coefficients for the appropriate $d$-variate monomials.
To identify near-resonances, we let $\mathbf{S}^{2:N}$ be the matrix
of integers whose columns are the powers of the $d$-variate monomials
from order 2 to $N$. We then define a matrix $\boldsymbol{\Delta}^{2:N}$
containing all relevant integer linear combinations of eigenvalues
as follows:
\begin{equation}
(\boldsymbol{\Delta}^{2:N})_{j,k}=(\mathrm{Im}\boldsymbol{\Lambda})_{j,j}-\sum_{s=1}^{d}(\mathrm{Im}\boldsymbol{\Lambda})_{s,s}(\mathbf{S}^{2:N})_{s,k}.\label{eq:innerresonances}
\end{equation}

Following the approach used in universal unfolding principles \cite{murdock2003},
we collect in a set $S$ the row and column indices of the entries
of $\boldsymbol{\Delta}^{2:N}$ for which near-resonances occur, i.e.,
for which the corresponding entry of $\boldsymbol{\Delta}^{2:N}$
is smaller in norm than a small, preselected threshold. (The default
threshold is $10^{-8}$ in \texttt{SSMLearn}.) The entries of $\mathbf{H}$
and $\mathbf{H}_{\star}$ with indices contained in $S$ are then
set to zero but the corresponding monomial terms are retained in $\mathbf{n}(\mathbf{z};\mathbf{N})$.
Conversely, coefficients of non-near-resonant entries of $\mathbf{H}$
and $\mathbf{H}_{\star}$ are selected in a way so that the corresponding
non--near-resonant monomials vanish from the normal form $\mathbf{n}(\mathbf{z};\mathbf{N})$.
As a result, the matrix $\mathbf{N}$ is sparse, containing only the
coefficients of essential, near-resonant monomials. 

For example, if $d=2$, $N=3$ and the eigenvalues of $D\mathbf{r}(\boldsymbol{0})$
form a complex pair $\lambda=\alpha_{0}\pm i\omega_{0}$ with $\omega_{0}=\mathcal{O}(1)$,
then we have
\begin{equation}
\mathbf{S}^{2:N}=\begin{bmatrix}2 & 1 & 0 & 3 & 2 & 1 & 0\\
0 & 1 & 2 & 0 & 1 & 2 & 3
\end{bmatrix},\,\,\boldsymbol{\Delta}^{2:N}=\begin{bmatrix}-\omega_{0} & \omega_{0} & 3\omega_{0} & -2\omega_{0} & 0 & 2\omega_{0} & 4\omega_{0}\\
-3\omega_{0} & -\omega_{0} & \omega_{0} & -4\omega_{0} & -2\omega_{0} & 0 & 2\omega_{0}
\end{bmatrix}.\label{eq:innerresonance2DSSMcubic}
\end{equation}
Only two elements of $\boldsymbol{\Delta}^{2:N}$ are (near-) zero,
and hence the reduced dynamics in extended normal form will require
learning the following coefficients:
\begin{equation}
\mathbf{H}_{\star}=\begin{bmatrix}h_{20} & h_{11} & h_{02} & h_{30} & 0 & h_{12} & h_{03}\\
\bar{h}_{02} & \bar{h}_{11} & \bar{h}_{20} & \bar{h}_{03} & \bar{h}_{12} & 0 & \bar{h}_{30}
\end{bmatrix},\,\,\,\,\mathbf{N}=\begin{bmatrix}0 & 0 & 0 & 0 & h_{21} & 0 & 0\\
0 & 0 & 0 & 0 & 0 & \bar{h}_{21} & 0
\end{bmatrix}.\label{eq:cubicnormalformcoefficients}
\end{equation}
The corresponding cubic polar form (\ref{eq:cubicnormalform}) is
then obtained from the relations $\mathbf{z}=(\rho e^{i\theta},\rho e^{-i\theta})$
and $h_{21}=\beta+i\gamma$.

For a data-driven construction of the extended normal form \eqref{eq:normalformSSMs},
we first obtain an estimate for the Jacobian $D\mathbf{r}(\boldsymbol{0})$
from linear regression. This determines the matrix $\mathbf{B}$ and
the types of monomials arising in $\mathbf{h}^{-1}$ and $\mathbf{n}$.
Next, we note that the flow map $\boldsymbol{\varphi}_{\mathbf{r}}^{t}$
of the SSM-reduced dynamics and the flow map $\boldsymbol{\varphi}_{\mathbf{n}}^{t}$
of the extended normal form are connected through the conjugacy relationship
$\boldsymbol{\varphi}_{\mathbf{n}}^{t}=\mathbf{h}^{-1}\circ\boldsymbol{\varphi}_{\mathbf{r}}^{t}\circ\mathbf{h}$.
We find the nonzero complex coefficients of $\mathbf{h}^{-1}$ and
$\mathbf{n}$ by minimizing the error in this exact conjugacy over
the available $P$ data points, represented in the $\mathbf{\boldsymbol{\eta}}$
coordinates. Specifically, we determine the nonzero elements of $\mathbf{H}_{\star}$
and $\mathbf{N}$ as 
\begin{equation}
\begin{split}(\mathbf{H}_{\star}^{*},\mathbf{N}^{*})= & {\displaystyle \mathrm{arg}\min_{\mathbf{H}_{\star},\mathbf{N}}\sum_{j=1}^{P}\left\Vert \frac{d}{dt}\mathbf{h}^{-1}(\boldsymbol{\eta}_{j};\mathbf{H}_{\star})-\mathbf{n}\left((\mathbf{h}^{-1}(\boldsymbol{\eta}_{j};\mathbf{H}_{\star});\mathbf{N}\right)\right\Vert ^{2},}\\
 & (\mathbf{N})_{s,k}=0,\,\,\forall(s,k)\in S;\,\,(\mathbf{H}_{\star})_{s,k}=0,\,\,\forall(s,k)\notin S.
\end{split}
\label{eq:learningdynamics}
\end{equation}
 Once $\mathbf{h}^{-1}$ is known, we obtain the coefficients $\mathbf{H}$
of $\mathbf{h}$ via regression.

As initial condition for the minimization problem \eqref{eq:learningdynamics},
we set all unknown coefficients to zero. This initial guess assumes
linear dynamics, which the minimization corrects as needed. We can
compute the time derivative in (\ref{eq:learningdynamics}) reliably
using finite differences, provided that the sampling time $\Delta t$
of the trajectory data is small compared to the fastest timescale
of the SSM dynamics. For larger sampling times, one should use the
discrete formulation of SSM theory, as discussed in section \ref{subsec:Existence-of-SSMs}
and \cite{haller16}. In that formulation, the conjugacy error must
be formulated for the 1-step prediction error of the normal form flow
map $\boldsymbol{\varphi}_{\mathbf{n}}^{\Delta t}(\mathbf{z})$. The
matrix defined in eq. (\ref{eq:innerresonances}) also carries over
to the discrete time setting, with $\boldsymbol{\Lambda}$ defined
as the diagonal matrix of the logarithms of the eigenvalues of $D\boldsymbol{\varphi}_{\mathbf{r}}^{\Delta t}(\boldsymbol{0})$. 

\subsection{Prediction of forced response from unforced training data\label{subsec:Prediction-of-forced-response-methods}}

Forced SSMs continue to be embedded in our observable space, provided
that we also include the phase of the forcing among our observables
\cite{stark99}. (In the simplest case of periodic forcing, this
inclusion is not necessary, as we pointed out Section \ref{subsec:Embedding-of-SSMs}).
The quasiperiodic SSM-reduced normal form of system \eqref{eq:1storder_system}
in the observable space takes the general form 
\begin{equation}
\begin{aligned}\dot{\rho}_{j} & =\alpha_{j}(\boldsymbol{\rho},\boldsymbol{\theta})\rho_{j}-\sum_{\mathbf{k}\in K_{j}^{\pm}}f_{j,\boldsymbol{\mathbf{k}}}\sin\left(\langle\mathbf{k},\boldsymbol{\Omega}\rangle t+\phi_{j,\mathbf{k}}\mp\theta_{j}\right),\\
\dot{\theta}_{j} & =\omega_{j}(\boldsymbol{\rho},\boldsymbol{\theta})+\sum_{\mathbf{k}\in K_{j}^{\pm}}\frac{f_{j,\mathbf{k}}}{\rho_{j}}\cos\left(\langle\mathbf{k},\boldsymbol{\Omega}\rangle t+\phi_{j,\mathbf{k}}\mp\theta_{j}\right),
\end{aligned}
\,\,\,\,\,\,\,j=1,2,...m,\,\,\,\,\,\,\mathbf{k}\in\mathbb{Z}^{\ell},\,\,\,\,\,\,\boldsymbol{\Omega}\in\mathbb{R}_{+}^{\ell},\label{eq:quasiperiodicSSMforcing}
\end{equation}
where the terms $f_{j,\mathbf{k}}$ and $\phi_{j,\mathbf{k}}$ are
the forcing amplitudes and phases for each mode of the SSM and for
each forcing harmonic $\langle\mathbf{k},\boldsymbol{\Omega}\rangle$,
while $K_{j}^{\pm}$ are the set containing the indexes $\mathbf{k}$
of the resonant forcing frequencies for mode $j$ (see the Supplementary
Information). The normal form \eqref{eq:quasiperiodicSSMforcing}
will capture non-linearizable dynamics arising from resonant interactions
between the eigenfrequencies of the spectral subspace $E$ (which
may also contain internal resonances) and the external forcing frequencies
in $\boldsymbol{\Omega}$. One can use numerical continuation \cite{dankowicz2013}
to find nontrivial co-existing steady states (such as periodic orbits
and invariant tori) in eq. \eqref{eq:quasiperiodicSSMforcing} under
varying forcing amplitudes and forcing frequencies. 

To predict forced response from the SSM-based model trained on unforced
data, the forcing amplitude $f$ relevant for eq. \eqref{eq:periodicSSMforcing}
in the observable space needs to be related to the forcing amplitude
$\left|\epsilon\mathbf{f}_{1}\right|$ relevant for system \eqref{eq:1storder_system}
in the physical phase space. This involves (1) employing a single
forcing amplitude-frequency pair $\left(\left|\epsilon\mathbf{f}_{1}\right|,\Omega\right)$
in the experiment (2) measuring the periodic observable response $\mathbf{y}(t)$
(3) computing the corresponding normalized reduced and normalized
response amplitude $\rho_{0}$ (4) substituting $\rho_{0}$ into the
first formula in \eqref{eq:closedformFRC} and (5) solving for $f$
in closed form. This $f$ can then be used to make a prediction for
the full FRC and response phase via \eqref{eq:closedformFRC} in the
experiment for arbitrary $\Omega$ forcing frequencies. The predicted
FRC may have several connected components, including isolated responses
(isolas) that are notoriously difficult to detect by numerical or
experimental continuation \cite{ponsioen2019}. 

\subsection{Summary of the algorithm\label{subsec:Summary-of-the-algorithm}}

The data-driven model reduction method used in this paper is available
in the open-source \textsc{MATLAB}\textsuperscript{®} package \texttt{SSMLearn}.
User input is the measured trajectory data of the autonomous dynamical
system ($\epsilon=0$), the SSM dimension $d$, the polynomial orders
or approximation $(M,N)$ for the SSM and for the extended normal
form, as well as the type of the dynamical system (discrete or continuous).
If the number of observables is not sufficient for manifold embedding,
the data is automatically augmented with delays to reach the minimum
embedding dimension $p=2d+1$. If the manifold learning returns poor
results (due to, e.g., insufficient closeness of the data to the SSM),
then the starting value of $p$ can be increased until a good embedding
is found. Then, the algorithm learns the SSM geometry in observable
space and, after unsupervised detection of the required normal form,
identifies the extended normal form of the reduced dynamics. The level
of accuracy can be increased with larger polynomial orders, keeping
in mind that excessive orders may lead to overfitting.

\texttt{SSMLearn} also offers all the tools we have used in this paper
to analyze the reduced dynamics and make predictions for forced response
from unforced training data. In particular, it contains the \textsc{MATLAB}\textsuperscript{®}-based
numerical continuation core \textsc{coco} \cite{dankowicz2013}.
which can compute steady state and help with the design of nonlinear
control strategies. In principle, there are no restrictions on the
dimensions of the reduced order model, yet the larger the SSM is,
the more computationally expensive the problem becomes. 

Qualitative or partial \textit{a priori} knowledge of the linearized
dynamics (e.g., some linearized modes and frequencies) helps in finding
good initial conditions for trajectories to be used in \texttt{SSMLearn}.
For example, the resonance decay method \cite{peeters2011a} (which
we exploited in our sloshing example), targets a specific 2-dimensional,
stable SSM in laboratory experiments. This method consists of empirically
isolating a resonant periodic motion on the SSM based on its locally
maximal amplitude response under a forcing frequency sweep. Discontinuing
the forcing will then generate transient decay towards the equilibrium
in a close proximity of the SSM. For noisy data, filtering or dimensionality
reduction can efficiently de-noise the data \cite{bishop2006}, provided
that the polynomial orders used for the description of the SSM and
its reduced dynamics are not excessively large (see the Supplementary
Information). For higher-dimensional SSMs, it is desirable to collect
diverse trajectories to avoid bias towards specific motions. Good
practice requires splitting the data sets into training, testing and
validation parts.
\begin{algorithm}
\caption{\texttt{SSMLearn}}

\textbf{Input parameters}: SSM dimension $d$, polynomial approximation
orders $(M,N)$, selection among discrete or continuous time dynamics

\textbf{Input data}: measured unforced trajectories

\textbf{Output}: SSM geometry, extended normal form of reduced dynamics,
predictions for forced response.

\textbf{1} Embed data in a suitable $p$-dimensional observable space
with $p>2d$.

\textbf{2} Identify the manifold parametrization in reduced coordinates.

\textbf{3} Estimate the normalized reduced dynamics after an automated
identification of the required type of extended normal form.

\textbf{4} Run analytics and prediction of forced response on the
SSM-reduced and normalized model.
\end{algorithm}

\subsubsection*{Acknowledgements}

B.B. and K.A. acknowledge financial support from the Deutsche Forschungsgemeinschaft (DFG, German Research Foundation) in the framework of the research unit FOR 2688 'Instabilities, Bifurcations and Migration in Pulsatile Flows' under Grant No. AV 156/1-1. K.A. acknowledges funding for an 'Independent Project for Postdocs' from the Central Research Development Fund of the University of Bremen.

\subsubsection*{Author contributions}

M.C. and G.H. designed the research. M.C. carried out the research.
M.C. and J.A. developed the software and analyzed the examples. B.B.
and K.A. performed the liquid sloshing experiments and participated
in their analysis. M.C. and G.H. wrote the paper. G.H. lead the research
team.

\subsubsection*{Competing interests}

The authors declare no competing interests.

\end{document}


\maketitle

\section{Details and additional analysis for the examples}

\subsection{Finite-element model of a damped-forced beam}

\subsubsection{Setup}

In our beam example \cite{cenedese21}, no body forces are present
and the straight equilibrium configuration is asymptotically stable.
The material properties of the beam are: the Young modulus is $70$
{[}GPa{]}, the density is $2700$ {[}kg/m$^{3}${]}, and the Poisson
ratio is $0.3$. The finite-element discretization is performed using
elements with cubic shape functions for the transverse deflection
and linear shape functions for the axial displacement \cite{jain2018}.
We use $12$ elements, which results in convergence in static and
dynamic simulations in the range of interest over uniform grid refinements.
The finite element model has $33$ degrees of freedom in total, including
transverse displacement, axial displacements and rotations. The slowest
eigenvalues of the linearized dynamics is $\lambda_{1}=-3.09\pm i657$,
giving rise to a two-dimensional slow spectral subspace $E_{1}$ with
spectral gap--the ratio between the real parts of the two slowest
eigenvalues--equal to $7$. Hence, the decay of faster linear modes
is more than seven times faster than that of the slowest mode.

\subsubsection{Comparison with \texttt{SSMTool} }

\begin{figure}[!t]
\includegraphics[width=1\textwidth]{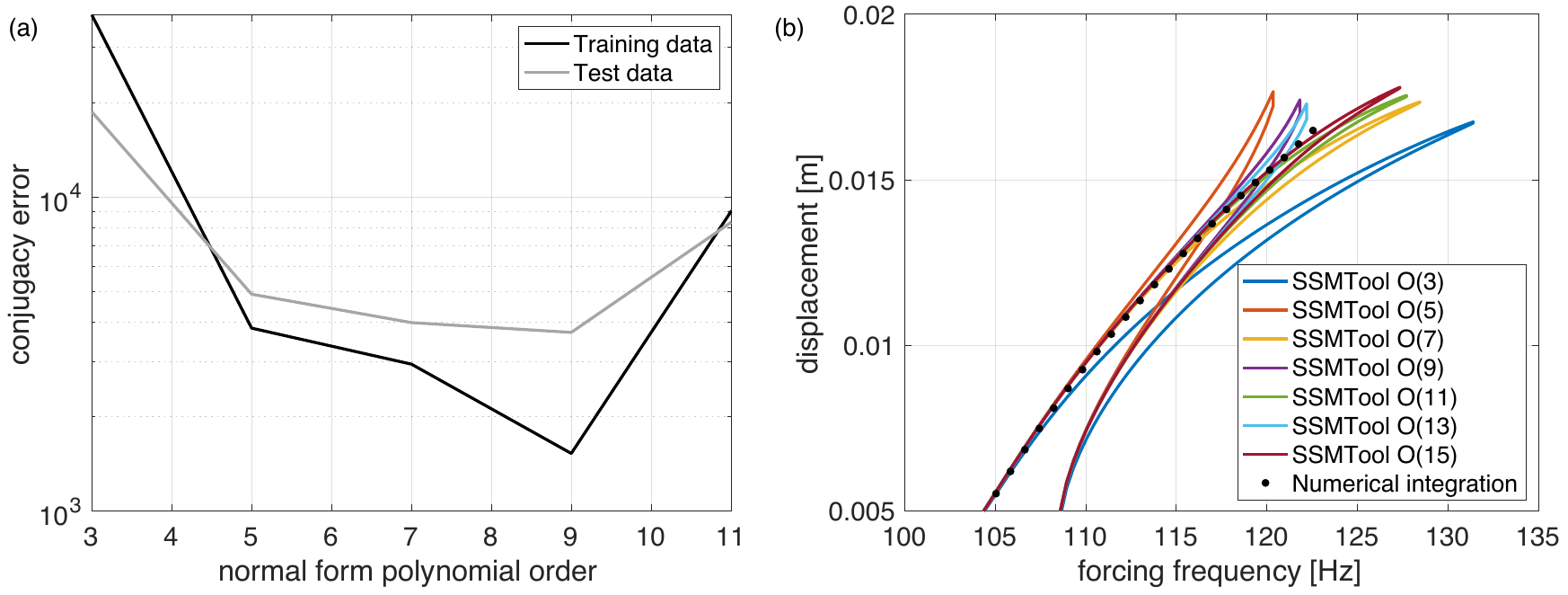}\caption{(a) Conjugacy error in the data-driven construction of the extended
normal form on the SSM of the damped-forced von Kármán beam as a function
of the polynomial order of the normal form. (b) FRCs obtained from
\texttt{SSMTool 2.0} of \cite{jain2021} for the beam at the highest
forcing amplitude analyzed in this paper, plotted for various polynomial
SSM-approximation orders. Also shown in black is the stable branch
of the FRC obtained from direct numerical integration.}\label{fig:vonkarmanbeamMore}
\end{figure}

We select the normal form order for the SSM-reduced model by minimizing
the conjugacy error defined in the \emph{Methods} section of \cite{cenedese21}.
This error is shown in Fig. \ref{fig:vonkarmanbeamMore}(a) as a function
of the polynomial order of approximation for the SSM, as computed
by \texttt{SSMLearn.} 

The training and testing trajectory data set comprises two trajectories,
with the training one featuring higher amplitudes than the testing
one. The conjugacy error reaches its global minimum for both trajectories
at the normal-form order $9$, as seen in Fig. \ref{fig:vonkarmanbeamMore}(a).
We nevertheless select the order 7, which produces slightly higher
errors but yields a simpler reduced model.

The slow convergence in Figure \ref{fig:vonkarmanbeamMore}(b) shows
that the results from the \texttt{SSMTool 2.0} of \cite{jain2021}
are approaching the boundary of the domain of convergence of the Taylor
expansion for the SSM. Indeed, even for an $\mathcal{O}(15)$ approximation,
the analytically predicted frequency responses do not yet match the
results from numerical integration, although the prediction error
decreases under increasing polynomial order. In contrast, as already
noted in our discussion of this example in \cite{cenedese21}, \texttt{SSMLearn}
returns more accurate results already at $\mathcal{O}(7)$ thanks
to its data-driven approach.

\subsubsection{Performance under noisy observable data}

\begin{figure}[!t]
\includegraphics[width=1\textwidth]{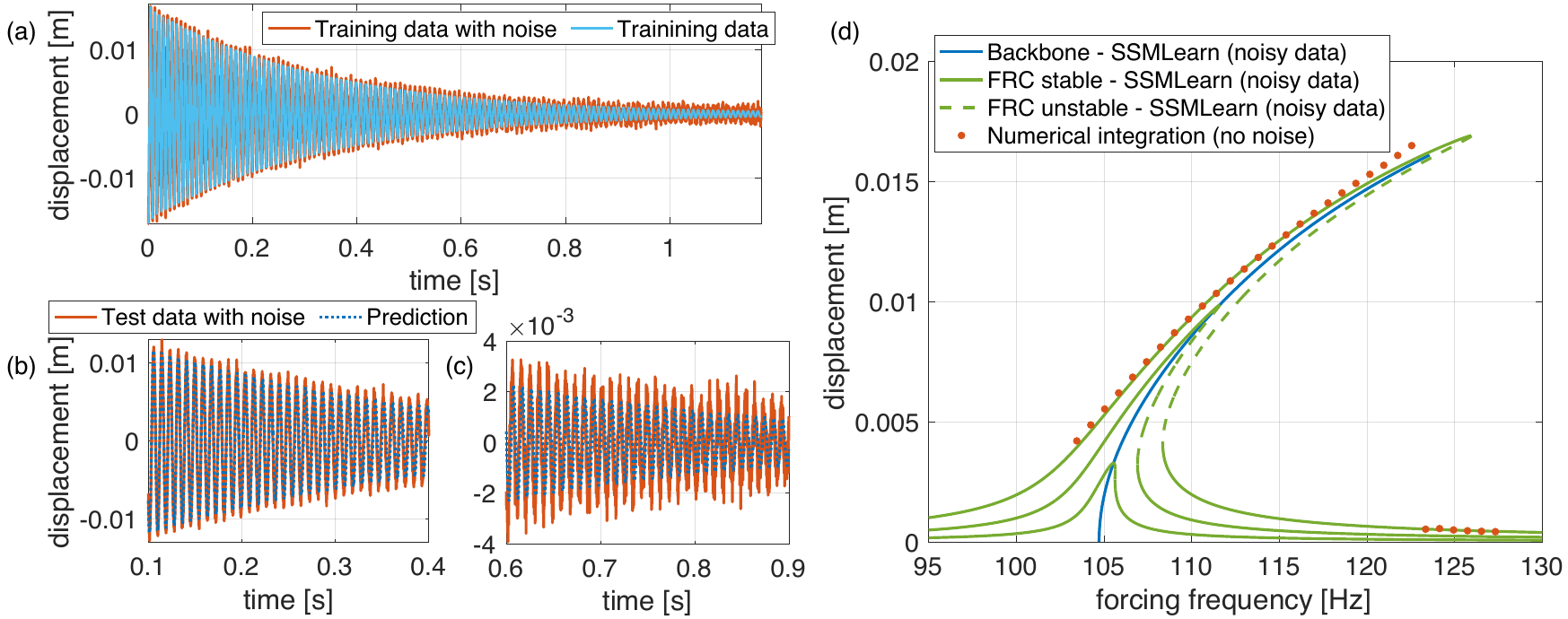}\caption{(a) Noisy training trajectory and its noise-free counterpart for the
von Kármán beam. (b) Noisy testing trajectory and its reconstruction
from the SSM-reduced model trained on noisy data. (c) Same as (b)
but zoomed in on the initial part of the time series. (d) FRCs predicted
from noisy observer data and compared to direct numerical simulations
in the absence of noise.}\label{fig:vonkarmanbeamnoisy}
\end{figure}

To illustrate the robustness of SSM-reduced models under noise, we
perturb the decaying training data of the damped beam with white noise
of variance 0.5 {[}mm{]}, as shown in Fig. \ref{fig:vonkarmanbeamnoisy}(a).
This emulates measurement noise or uncertainty for the observable.

The available theory of SSMs also holds approximately for finite times
in this noisy setting, given that any realization of a noisy dynamical
system can be approximated arbitrarily closely by a quasiperiodic
dynamical system over a fixed finite time interval. The noisy attractor
born out of the equilibrium under forcing can therefore be approximated
by a quasiperiodic attractor with sufficiently many frequencies. Due
to this multi-frequency time dependence of the corresponding SSM,
we need to increase the dimensions of the delay-embedding space used
for the noisy data. Accordingly, we choose to set the embedding dimension
to be 200 (the scalar observable plus 199 delays) and continue to
use the flat manifold approximation with an $\mathcal{O}(7)$ reduced
dynamics, as we did in the case of perfect measurements. As seen in
Fig. \ref{fig:vonkarmanbeamnoisy}(b), this approach practically filters
out noise from our reduced model coordinates. Indeed, the SSM-reduced
model built from noisy data still reconstructs trajectories with an
error of only $\mathrm{NMTE}=7.6$ \%, as shown in Figs. \ref{fig:vonkarmanbeamnoisy}(b,c).
Most importantly, the FRCs predicted from noisy data still align with
the noise-free FRC obtained from direct numerical integration, as
shown in Fig. \ref{fig:vonkarmanbeamnoisy}(d).

\subsection{Vortex-shedding behind a cylinder}

\subsubsection{Setup}

For the vortex-shedding example in \cite{cenedese21}, we set the
fluid kinematic viscosity to 0.01 {[}m$^{2}$/s{]} and the inflow
horizontal velocity to 0.7 {[}m/s{]}, resulting in a Reynolds number
equal to 70. We simulate the flow using the Python-based computational
platform FEniCS \cite{FenicsBook2012}. The mesh is formed by triangular
elements of the Lagrange family; the resulting discretized model has
a phase space $\mathbb{R}^{n}$ with dimension $n=76876$. We integrate
the flow in time using a modified version of Chorin's method \cite{FenicsTutorial2016}
with time step 0.02 {[}s{]}. The linear stability analysis of the
steady solution is performed via Krylov--Arnoldi iterations \cite{loiseau19},
while the POD modes (see \cite{holmes12}) of the limit cycle are
computed using the snapshot method \cite{taira2107} applied to the
velocities and pressures along the limit cycle. For the vortex shedding,
the two leading POD modes capture most of the energy content of the
flow on the limit cycle, as found originally by \cite{noack2003}.

We use nine trajectories for training the SSM-reduced model, with
initial conditions that are small perturbations of the steady solution
along its unstable directions. Specifically, we set
\begin{equation}
\mathbf{x}(0)=a\mathbf{w}_{u,1}\cos\beta+a\mathbf{w}_{u,2}\sin\beta,
\end{equation}
where $\mathbf{w}_{u,1},\,\mathbf{w}_{u,2}$ are real unit vectors
spanning the unstable subspace, $a=0.2$ and the angle $\beta$ attains
9 uniformly spaced values in the interval $[0,2\pi)$.

\subsubsection{Strong deformation along the SSM}

\begin{figure}[!t]
\includegraphics[width=1\textwidth]{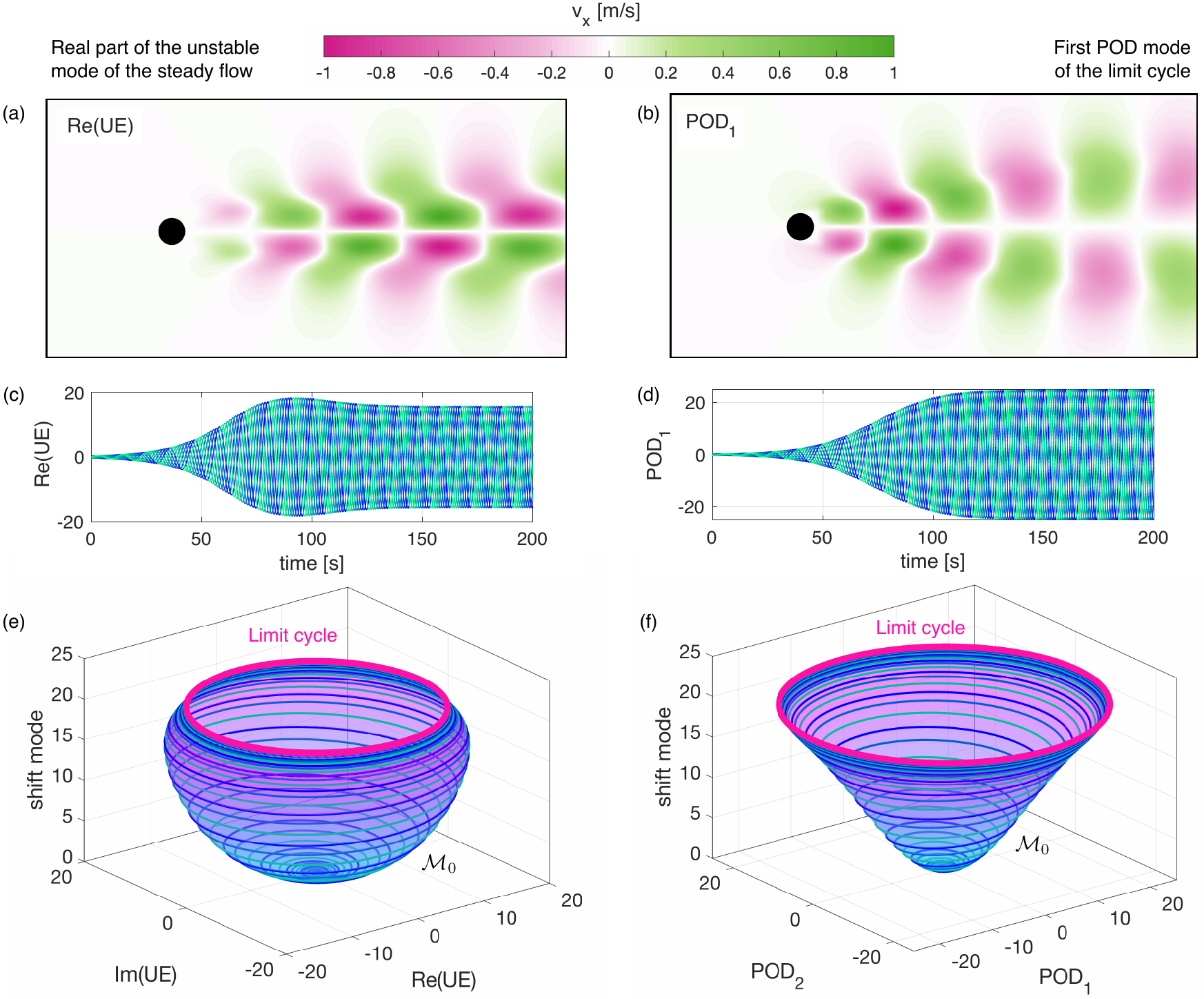}\caption{Modes and coordinates for the dynamics along the SSM (unstable manifold)
of the vortex shedding example. (a) Deviation of the horizontal velocity
from the steady flow for the real part of the unstable mode of the
unstable eigenspace. (b) Same for the first POD mode of the limit
cycle. (c,d) Projection of the full trajectory onto the modes shown
in (a) and (b), respectively. (e) The SSM (with some reduced trajectories)
plotted over the unstable eigenspace UE and the shift mode, as defined
by \cite{noack2003}. (f) The same SSM plotted over the two-dimensional
leading POD subspace and the shift mode.}\label{fig:vortexsheddingModes}
\end{figure}

As shown in Fig. 5(f) in \cite{cenedese21}, the SSM
(unstable manifold) develops a fold over the unstable subspace of
the fixed point representing the steady wake flow. Figure \ref{fig:vortexsheddingModes}
shows additional details of this phenomenon.
Trajectories initialized close to the unstable fixed
point initially follow its unstable eigenspace (UE), whose real part
is shown in Fig. \ref{fig:vortexsheddingModes}(a). Continuing to
evolve along the SSM, these trajectories converge to a limit cycle,
whose first POD mode is shown in Fig. \ref{fig:vortexsheddingModes}(b).
As seen from a bump in the envelope of the signal in Fig. \ref{fig:vortexsheddingModes}(c),
the SSM develops a fold over UE. At the same time, there is no similar
bump over the leading POD modes in Fig. \ref{fig:vortexsheddingModes}(d),
which suggest a lack of a fold over the POD modes.

The actual geometry of the fold over the UE subspace
is shown in Fig. \ref{fig:vortexsheddingModes}(e), with the SSM
plotted in the space of three coordinates: the real and imaginary
parts of the corresponding complex unstable eigenvector of the linearization,
and the shift mode defined by \cite{noack2003}. In contrast, no
such fold is present in Fig. \ref{fig:vortexsheddingModes}(f), in
which the same SSM is plotted using the two leading POD modes and
the shift mode as coordinate axes.

\subsubsection{Comparison with SINDy and DMD}

\begin{figure}[!t]
\includegraphics[width=1\textwidth]{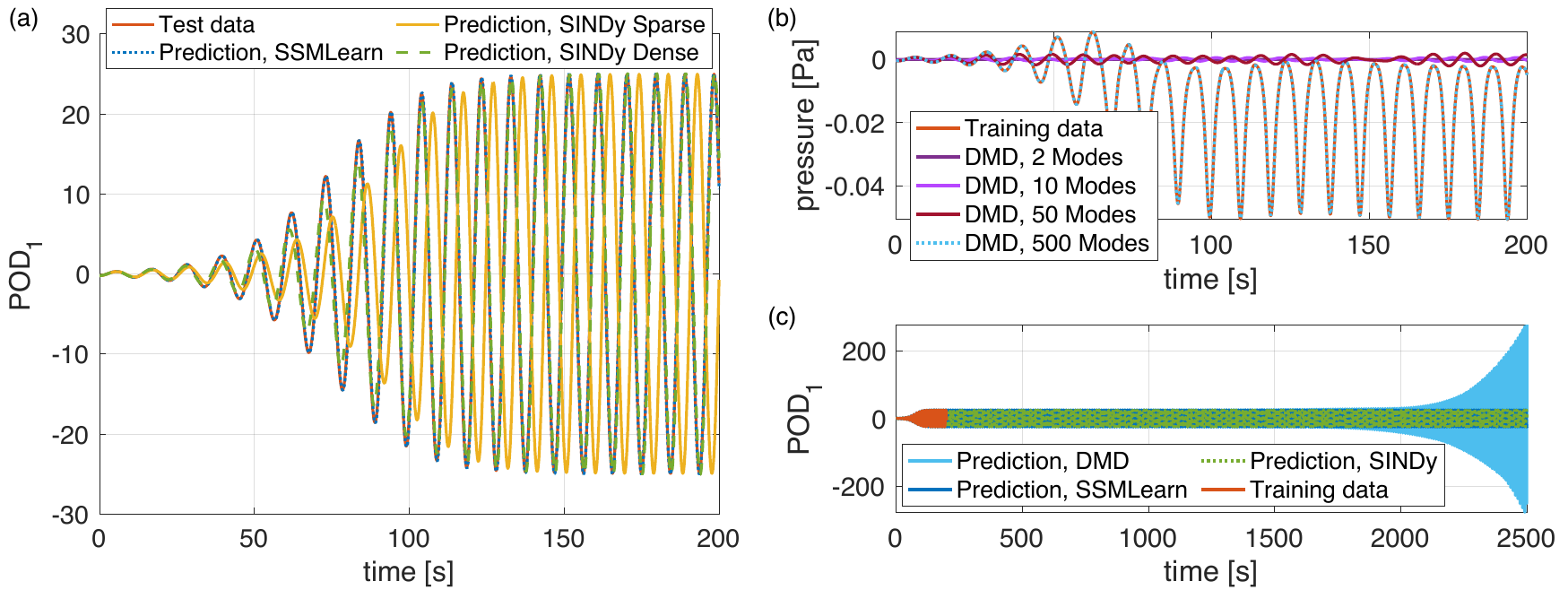}\caption{Comparison of \texttt{SSMLearn} with two popular data-driven model
identification methods on the vortex shedding problem. (a) Comparison
of the predictions of\texttt{ SSMLearn} and SINDy for the projection
of the model testing trajectory to the first POD mode. (b) Reconstruction
of the pressure history along the training data (see Fig. 5e of the
paper) via DMD using various numbers of DMD modes. (c) Predictions
from the three reduced-order models (\texttt{SSMLearn} at order 11,
SINDy dense, DMD with 500 modes) for the training trajectory over
a longer time interval that exceeds the training window.}\label{fig:vortexsheddingSINDyDMD}
\end{figure}

Figure \ref{fig:vortexsheddingSINDyDMD} compares results produced
by \texttt{SSMLearn} with those from two leading data-driven model
identification methods applied to the same dataset. 

The first of these methods, SINDy \cite{brunton2016}, uses sparse
regression to approximate an envisioned reduced vector field representation
of a nonlinear dynamical system. As this method is not originally
intended for model reduction, we examine its performance after we
have projected the eight training trajectories to the two leading
POD modes computed along the limit cycle of this system.

When used with polynomial regression, SINDy has two input parameters:
the polynomial order and a regularizing parameter $\lambda$ that
controls the sparsity in the vector field to be learned. For the choice
$\lambda=3.16\cdot10^{-5}$ , SINDy returns the sparse polynomial
model 
\begin{equation}
\begin{array}{rl}
\dot{\eta}_{1}= & 0.00367+0.0540\eta_{1}-0.545\eta_{2}-10^{-5}\eta_{1}(8.77\eta_{1}^{2}+8.45\eta_{2}^{2})-10^{-4}\eta_{2}(1.93\eta_{1}^{2}+1.92\eta_{2}^{2}),\\
\dot{\eta}_{2}= & 0.00175+0.534\eta_{1}+0.0394\eta_{2}-10^{-5}\eta_{2}(6.32\eta_{1}^{2}+6.33\eta_{2}^{2})+10^{-4}\eta_{1}(2.22\eta_{1}^{2}+2.19\eta_{2}^{2}),
\end{array}\label{eq:vortexsheddingSINDysparse}
\end{equation}
where $(\eta_{1},\eta_{2})$ are coordinates along the two leading
POD modes of the limit cycle. Equation (\ref{eq:vortexsheddingSINDysparse})
has a sizable reconstruction error of $\mathrm{NMTE}=108\%$, which
is clearly visible in Fig. \ref{fig:vortexsheddingSINDyDMD}(a). \cite{fukami21}
obtain comparably large errors for the same problem after applying
SINDy to a latent (i.e., reduced) representation of the dynamics generated
by an auto-encoder.

Cross-validation yields the best hyper-parameter value $\lambda=1.08\cdot10^{-8}$,
which reduces the reconstruction error to $\mathrm{NMTE}=10.2\,\%$.
The resulting model returned by SINDy, however, is dense: all its
polynomial coefficients (up to degree 5) are nonzero. This is to be
contrasted with the reconstruction error $\mathrm{NMTE}=3.86$ $\%$
of the sparse SSM-reduced model of \cite{cenedese21}. The simultaneous
sparsity and accuracy of the SSM-reduced model illustrates the benefits
of exploiting the precise knowledge of the phase space geometry of
the underlying dynamical system.

As a second method to compare with\texttt{ SSMLearn}, we select the
DMD \cite{schmid10,kutz16}. Because this algorithm fundamentally
assumes near-linear dynamics, using several training trajectories
with their different levels of nonlinearity leads to poor reduced
models. To improve the performance of DMD, we therefore use a single
training trajectory with sizable displacements. Figure \ref{fig:vortexsheddingSINDyDMD}(b)
shows that while low-rank DMD fits fail to predict the trajectory
from the initial condition, the DMD fit becomes very accurate for
500 modes. This fundamentally linear DMD model, however, cannot capture
the limit cycle, as evidenced by the ultimate convergence of the DMD
prediction to infinity. This is seen in Fig. \ref{fig:vortexsheddingSINDyDMD}(c),
which shows the predictions for the testing trajectory from all three
approaches over a longer time interval.

\subsection{Fluid sloshing experiment}

Higher-amplitude forcing not considered in \cite{cenedese21} can
lead to an apparent $1:3$ resonance in the forced response of the
fluid surface \cite{bauerlein2021}. \texttt{SSMLearn} is fully equipped
to handle such a resonant interaction, but the available experimental
data does not provide sufficient information about the higher modes
involved in the interaction. More focused resonance decay experiments
exciting those higher modes at higher amplitudes are needed for \texttt{SSMLearn
}to construct an accurate reduced model on a four-dimensional, resonant
SSM. This is the subject of our ongoing work, to be reported elsewhere.

\begin{figure}[!t]
\includegraphics[width=1\textwidth]{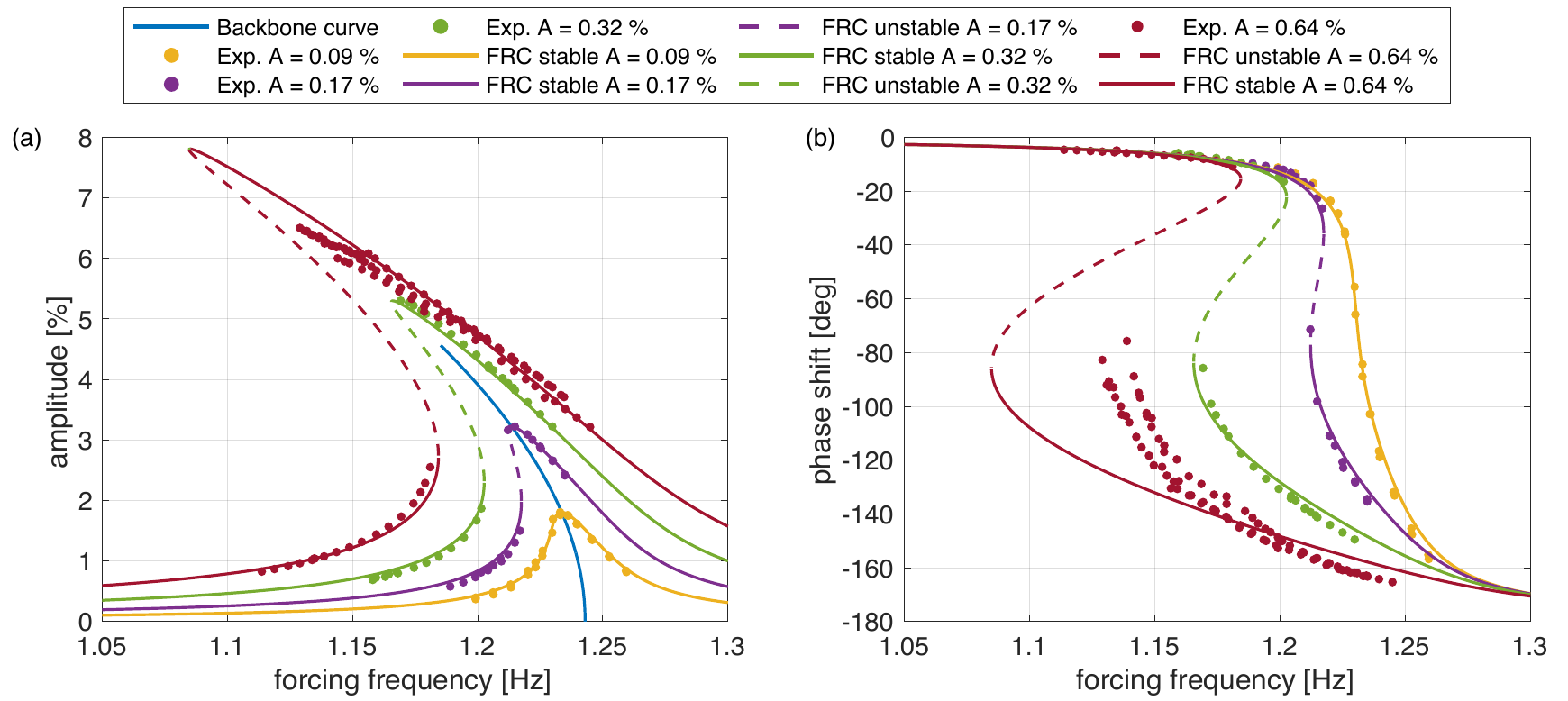}\caption{Inclusion of period-three resonant response (red dots) observed from
the fluid sloshing experiment for amplitude (a) and phase (b). Red
lines indicate the predicted response from the two-dimensional SSM
calculation at the same forcing level.}\label{fig:sloshing with resonance}
\end{figure}

The four-dimensional SSM-reduced model is expected to comprise a second
degree of freedom that is weakly coupled to our current two-dimensional
SSM model (eq. (12) in \cite{cenedese21}). Indeed, the two-dimensional
SSM model in \cite{cenedese21} already captures the resonant branch
of the FRC curve family with surprisingly high overall accuracy (see
Fig. \ref{fig:sloshing with resonance}). This one-degree-of-freedom
model, however, cannot account for the slight periodic variation in
the experimentally measured response amplitudes, manifested by the
three recurrent, distinct vertical dots in Fig. \ref{fig:sloshing with resonance}
for each forcing frequency in the resonant domain.

\subsubsection{Forced response via DMD}

\begin{figure}[!t]
\includegraphics[width=1\textwidth]{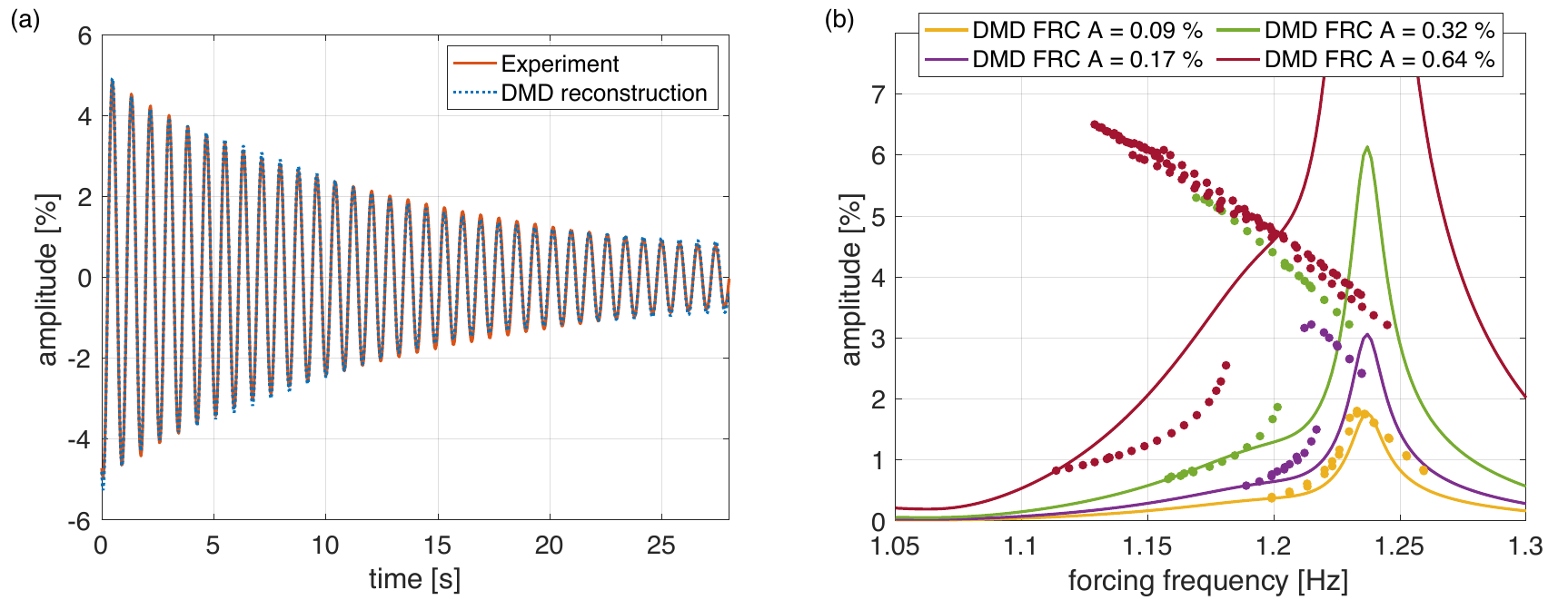}\caption{Attempted prediction of the forced sloshing response from a DMD model
fitted to unforced data. (a) With a high enough delay embedding dimension,
DMD accurately matches the decaying center-of-mass signal of the unforced
fluid mass. (b) Adding calibrated forcing to the DMD model, however,
will always yield linear response. While approximately correct for
very low forcing amplitudes, the forced DMD prediction therefore necessarily
fails at higher amplitudes in the non-linearizable regime.}\label{fig:sloshingDMD}
\end{figure}

To demonstrate that sloshing, as a non-linearizable problem, lies
outside the applicability of liner data-driven modeling, we also attempt
to predict the forced response using DMD \cite{schmid10,kutz16}.
To help DMD discover different modal contents arising from nonlinearities,
we choose a longer delay of 10 time steps with a delay-embedding dimension
of 20 to create a characteristically non-flat SSM in the embedding
space (see Section \ref{sec:Flat SSM in delay embedding} for related
information). This choice results in a close match of the decaying
training signal with four DMD modes as shown in Fig. \ref{fig:sloshingDMD}(a).
This multimodal linear DMD model approximates the nonlinear signal
as a superposition of high- and low-amplitude trajectory patterns
(see, e.g., \cite{dylewsky21}).

To obtain a forced-response prediction from this DMD model, we simply
add an appropriate forcing vector to the discrete, linear, unforced
DMD model and iterate the resulting inhomogeneous linear difference
equation to obtain its steady-state solution. We obtain the forcing
vector used for this calculation by projecting a delay-embedded harmonic
load signal onto the four DMD modes. The magnitude of the modal forcing
is calibrated to the known experimental forced response curves at
the lowest amplitude point in the frequency sweep. The predicted periodic
response amplitudes obtained from the 4-mode forced DMD model are
shown in Fig. \ref{fig:sloshingDMD}(b) for the same forcing amplitudes
and frequencies used in Fig. \ref{fig:sloshing with resonance}(a).
The phase response predicted by DMD can also be computed using different
approximations and assumptions but yield widely differing results,
which we do not show here.

For very small forcing amplitudes, the sloshing dynamics is close
to linear. Accordingly, the prediction from the 4-mode forced DMD
model is adequate, as illustrated by the lowest amplitude response
curve in Fig. \ref{fig:sloshingDMD}(b). The forced response predicted
by the forced DMD model, however, completely misses the softening
nonlinearity of the forced response curve even for moderate forcing
amplitudes, let alone for large ones. This is no surprise, given that
the overhangs in the experimental forced response curves signal the
coexistence of multiple isolated stationary states, i.e., non-linearizable
dynamics.

\section{Geometry of SSMs in delay-embedding spaces\label{sec:Flat SSM in delay embedding}}

Both the first and the third example in \cite{cenedese21} involves
single scalar observables, and hence the embedding space for the SSMs
in these examples need to be constructed from delayed measurements
(see the \emph{Methods} of \cite{cenedese21}). As seen in Figs.
4 and 6 of \cite{cenedese21}, the SSM in the delay-embedding space
is close to a plane in both cases. As we will show next, this flatness
of the delay-embedding of the SSM turns out to be a general phenomenon
if a moderate number of small enough delays are employed in the embedding
of a signal that does not change too abruptly.

Specifically, for any smooth scalar time series, $s(t)$, and any
integer $k$, we can write 
\begin{equation}
s(t+k\Delta t)=s(t)+\dot{s}(t)k\Delta t+\mathcal{O}(\Delta t^{2}).
\end{equation}
As a consequence, a $p$-dimensional delay embedding vector, $\mathbf{y}\in\mathbb{R}^{p}$,
constructed from $s(t)$ can be written as 
\begin{equation}
\mathbf{y}(t)=\left(\begin{array}{c}
s(t)\\
s(t+\Delta t)\\
s(t+2\Delta t)\\
\vdots\\
s(t+(p-1)\Delta t)
\end{array}\right)=s(t)\left(\begin{array}{c}
1\\
1\\
1\\
\vdots\\
1
\end{array}\right)+\dot{s}(t)\left(\begin{array}{c}
0\\
1\\
2\\
\vdots\\
p-1
\end{array}\right)\Delta t+\mathcal{O}\left((p-1)^{2}\Delta t^{2}\right).\label{eq:delayembeddingexplained}
\end{equation}
If the time step $\Delta t$ and the number of delays $p$ is small
enough in \eqref{eq:delayembeddingexplained}, then $\mathbf{y}(t)$
always lies approximately on a plane spanned by two constant vectors.
In that case, therefore, a flat approximation to the SSM will always
be sufficient, as long as the higher-order derivatives of $s(t)$
are not too large, i.e., the signal does not change too abruptly relative
to the sampling time scale $\Delta t$. This general conclusion, however,
does not imply that the reduced dynamics on the nearly flat SSMs are
linear or that the SSM is nearly flat in the phase space of the underlying
dynamical system.

\section{The quasiperiodically forced, extended SSM normal form}

We recall the main subject of \cite{cenedese21}, an $n$-dimensional
dynamical system of the form 
\begin{equation}
\dot{\mathbf{x}}=\mathbf{A}\mathbf{x}+\mathbf{f}_{0}(\mathbf{x})+\epsilon\mathbf{f}_{1}(\mathbf{x},\boldsymbol{\Omega}t;\epsilon),\qquad\mathbf{f}_{0}(x)=\mathcal{O}(\left|x\right|^{2}),\qquad0\leq\epsilon\ll1,\label{eq:1storder_system}
\end{equation}
with a constant, semisimple matrix $\mathbf{A}\in\mathbb{R}^{n\times n},$
and with the class $C^{r}$ functions $\mathbf{f}_{0}\colon\mathcal{U}\to\mathbb{R}^{n}$
and $\mathbf{f}_{1}\colon\mathcal{U}\times\mathbb{T}^{\ell}\to\mathbb{R}^{n}$,
where $\mathbb{T}^{\ell}=S^{1}\times\ldots\times S^{1}$ is the $\ell$-dimensional
torus and $\boldsymbol{\Omega}\in\mathbb{R}_{+}^{l}$. The assumed
degree of smoothness for the right-hand side of (\ref{eq:1storder_system})
is $r\in\mathbb{N}^{+}\cup\left\{ \infty,a\right\} $ \cite{cabre03}.

We define the matrix $\mathbf{T}\in\mathbb{C}^{n\times n}$ whose
columns are the eigenvectors of $\mathbf{A}$. We also let $\boldsymbol{\Lambda}\in\mathbb{C}^{n\times n}$
denote the diagonal matrix of the corresponding eigenvalues, ordered
with decreasing real parts. We write the Fourier expansion of the
quasiperiodic forcing $\mathbf{f}_{1}(\mathbf{x},\boldsymbol{\Omega}t;\epsilon)$
in the form 
\begin{equation}
\begin{array}{cc}
\mathbf{f}_{1}(\mathbf{x},\boldsymbol{\Omega}t;\epsilon)={\displaystyle \sum_{\mathbf{k}\in\mathbb{Z}^{l}}\mathrm{Re}\left(\mathbf{f}_{\mathbf{k}}^{1}e^{i\langle\mathbf{k},\boldsymbol{\Omega}\rangle t}\right)+\mathcal{O}(\epsilon\|\mathbf{x}\|)}=\mathbf{T}\mathbf{g}_{1}(\boldsymbol{\Omega}t)+\mathcal{O}(\epsilon\|\mathbf{x}\|), & \mathbf{f}_{\mathbf{k}}^{1}\in\mathbb{C}^{n},\end{array}\label{eq:forcing_definition}
\end{equation}
with $\mathbf{g}_{1}(\boldsymbol{\Omega}t):=\mathbf{T}^{-1}\mathbf{f}_{1}(\mathbf{x},\boldsymbol{\Omega}t;0)$.
Using the modal coordinates $\mathbf{q}\in\mathbb{C}^{n\times n}$
defined by the relation $\mathbf{x}=\mathbf{T}\mathbf{q}$ , we rewrite
eq. (\ref{eq:1storder_system}) as the autonomous system 
\begin{equation}
\begin{array}{l}
\dot{\mathbf{q}}=\boldsymbol{\Lambda}\mathbf{q}+\mathbf{g}_{0}(\mathbf{q})+\epsilon\mathbf{g}_{1}(\boldsymbol{\varphi})+\mathcal{O}(\epsilon\|\mathbf{q}\|),\\
\dot{\boldsymbol{\varphi}}=\boldsymbol{\Omega},
\end{array}\label{eq:1storder_system_modal}
\end{equation}
where $\mathbf{g}_{0}(\mathbf{q})=\mathbf{T}^{-1}\mathbf{f}_{0}(\mathbf{Tq})$.

Focusing on the slow SSM of dimension $2m\leq n$ related to the first
$m$ complex conjugate pairs of eigenvalues (either all stable or
all unstable), we can write the extended normal form of the vector
field (\ref{eq:1storder_system_modal}) in the general polar form
\begin{equation}
\begin{aligned}\dot{\rho}_{j} & =\alpha_{j}(\boldsymbol{\rho},\boldsymbol{\theta})\rho_{j}-\sum_{\mathbf{k}\in K_{j}^{+}}f_{j,\mathbf{k}}\sin(\langle\mathbf{k},\boldsymbol{\Omega}\rangle t+\phi_{j,\mathbf{k}}-\theta_{j})-\sum_{\mathbf{k}\in K_{j}^{-}}f_{j,\mathbf{k}}\sin(\langle\mathbf{k},\boldsymbol{\Omega}\rangle t+\phi_{j,\mathbf{k}}+\theta_{j}),\\
\dot{\theta}_{j} & =\omega_{j}(\boldsymbol{\rho},\boldsymbol{\theta})+\sum_{\mathbf{k}\in K_{j}^{+}}\frac{f_{j,\mathbf{k}}}{\rho_{j}}\cos(\langle\mathbf{k},\boldsymbol{\Omega}\rangle t+\phi_{j,\mathbf{k}}-\theta_{j})+\sum_{\mathbf{k}\in K_{j}^{-}}\frac{f_{j,\mathbf{k}}}{\rho_{j}}\cos(\langle\mathbf{k},\boldsymbol{\Omega}\rangle t+\phi_{j,\mathbf{k}}+\theta_{j}),
\end{aligned}
\label{eq:quasiperiodicSSMforcing}
\end{equation}
for $j=1,2,...m$, where $K_{j}^{+}$ and $K_{j}^{-}$ are the sets
of resonant forcing frequencies for mode $j$ defined as 
\begin{equation}
\begin{array}{c}
K_{j}^{+}:=\left\{ {\displaystyle \mathbf{k}\in\mathbb{Z}^{l},\left|\lim_{\|\boldsymbol{\rho}\|\rightarrow0}\omega_{j}(\boldsymbol{\rho},\boldsymbol{\theta})-\langle\mathbf{k},\boldsymbol{\Omega}\rangle\right|\leq\delta}\right\} ,\\
K_{j}^{-}:=\left\{ {\displaystyle \mathbf{k}\in\mathbb{Z}^{l},\left|\lim_{\|\boldsymbol{\rho}\|\rightarrow0}\omega_{j}(\boldsymbol{\rho},\boldsymbol{\theta})+\langle\mathbf{k},\boldsymbol{\Omega}\rangle\right|\leq\delta}\right\} ,
\end{array}\label{eq:resonant_index}
\end{equation}
with $\delta$ being a small tolerance. Similar normal forms in prior
work \cite{szalai02,ponsioen2018,breunung2018,ponsioen2019,ponsioen2020,jain2021}
focus on periodic forcing in two-dimensional SSMs.

The detailed form of the expressions $\alpha_{j}(\boldsymbol{\rho},\boldsymbol{\theta})$
and $\omega_{j}(\boldsymbol{\rho},\boldsymbol{\theta})$, and the
Taylor expression of the SSM carrying the normal form \eqref{eq:quasiperiodicSSMforcing}
can be determined simultaneously via a recursive solution of the partial
differential equation emerging from the invariance of the SSM. To
obtain this invariance equation, we first consider a yet unknown parameterization
of the SSM in the form $\mathbf{q}=\mathbf{w}(\mathbf{z},\boldsymbol{\varphi};\epsilon)$
with $\mathbf{z}\in\mathbb{C}^{2m}$. The SSM-reduced dynamics will
then be of the general normal form $\dot{\mathbf{z}}=\mathbf{n}(\mathbf{z},\boldsymbol{\varphi};\epsilon)$,
satisfying $\mathbf{w}(\mathbf{0},\boldsymbol{\varphi};0)=\mathbf{0}$
and $\mathbf{n}(\mathbf{0},\boldsymbol{\varphi};0)=\mathbf{0}$.

If we substitute the parametrization and the reduced dynamics expressions
into eq. (\ref{eq:1storder_system_modal}), we obtain the invariance
equation \cite{haro16} 
\begin{equation}
D_{\mathbf{z}}\mathbf{w}(\mathbf{z},\boldsymbol{\varphi};\epsilon)\mathbf{n}(\mathbf{z},\boldsymbol{\varphi};\epsilon)+D_{\mathbf{\boldsymbol{\varphi}}}\mathbf{w}(\mathbf{z},\boldsymbol{\varphi};\epsilon)\boldsymbol{\Omega}=\boldsymbol{\Lambda}\mathbf{w}(\mathbf{z},\boldsymbol{\varphi};\epsilon)+\mathbf{g}_{0}(\mathbf{w}(\mathbf{z},\boldsymbol{\varphi};\epsilon))+\epsilon\mathbf{g}_{1}(\boldsymbol{\varphi})+\mathcal{O}(\epsilon\|\mathbf{w}\|).\label{eq:invariance_equation}
\end{equation}
By the smooth dependence of the SSM on $\epsilon$ \cite{ponsioen2020},
we can write $\mathbf{w}(\mathbf{z},\boldsymbol{\varphi};\epsilon)=\mathbf{w}_{0}(\mathbf{z})+\epsilon\mathbf{w}_{1}(\mathbf{z},\boldsymbol{\varphi})+\mathcal{O}(\epsilon^{2})$
and $\mathbf{n}(\mathbf{z},\boldsymbol{\varphi};\epsilon)=\mathbf{n}_{0}(\mathbf{z})+\epsilon\mathbf{n}_{1}(\mathbf{z},\boldsymbol{\varphi})+\mathcal{O}(\epsilon^{2})$.
Substituting these expression in eq. (\ref{eq:invariance_equation})
and the Taylor expansion $\mathbf{g}_{0}(\mathbf{w}(\mathbf{z},\boldsymbol{\varphi};\epsilon))=\mathbf{g}_{0}(\mathbf{w}_{0}(\mathbf{z}))+\epsilon D\mathbf{g}_{0}(\mathbf{w}_{0}(\mathbf{z}))\mathbf{w}_{1}(\mathbf{z},\boldsymbol{\varphi})+\mathcal{O}(\epsilon^{2})$,
we then equate terms of equal order in $\epsilon$ to obtain the invariance
equation to be detailed next.

\subsection{Solving the autonomous invariance equation}

Equating terms of order $\mathcal{O}(1)$ in eq. \eqref{eq:invariance_equation},
we obtain the autonomous invariance equation 
\begin{equation}
D\mathbf{w}_{0}(\mathbf{z})\mathbf{n}_{0}(\mathbf{z})=\boldsymbol{\Lambda}\mathbf{w}_{0}(\mathbf{z})+\mathbf{g}_{0}(\mathbf{w}_{0}(\mathbf{z})).\label{eq:invariance_equation_aut}
\end{equation}
We seek the solution of this partial differential equation via the
Taylor expansions 
\begin{equation}
\begin{array}{ccc}
\mathbf{w}_{0}(\mathbf{z})={\displaystyle \sum_{j=1}^{M}\mathbf{W}_{j}^{0}\mathbf{z}^{j}}, & \mathbf{n}_{0}(\mathbf{z})={\displaystyle \sum_{j=1}^{M}\mathbf{N}_{j}^{0}\mathbf{z}^{j}}, & \mathbf{g}_{0}(\mathbf{q})={\displaystyle \sum_{j=2}^{M_{g}}\mathbf{G}_{j}^{0}\mathbf{q}^{j}},\end{array}\label{eq:expansion_aut}
\end{equation}
where $\mathbf{z}^{j}$ denotes the family of all monomials in $2m$
variables of degree $j$, and $\mathbf{W}_{j}^{0},\mathbf{N}_{j}^{0}$
and $\mathbf{G}_{j}^{0}$ are matrices storing the corresponding polynomial
coefficients. We also introduce the notation $\mathbf{S}_{2m}^{j}$
for the matrix of integers whose columns are the exponents of the
$2m$-variate monomials of order $j$, i.e., the components of $\mathbf{z}^{j}$.
For example, if $m=2$ and $j=3$ , then 
\begin{equation}
\begin{array}{c}
\mathbf{z}^{2}=(z_{1}^{2},z_{1}z_{2},z_{1}z_{3},z_{1}z_{4},z_{2}^{2},z_{2}z_{3},z_{2}z_{4},z_{3}^{2},z_{3}z_{4},z_{4}^{2})^{T},\\
\mathbf{S}_{2m}^{2}=\left[\begin{array}{cccccccccc}
2 & 1 & 1 & 1 & 0 & 0 & 0 & 0 & 0 & 0\\
0 & 1 & 0 & 0 & 2 & 1 & 1 & 0 & 0 & 0\\
0 & 0 & 1 & 0 & 0 & 1 & 0 & 2 & 1 & 0\\
0 & 0 & 0 & 1 & 0 & 0 & 1 & 0 & 1 & 2
\end{array}\right].
\end{array}
\end{equation}

After substituting the expression (\ref{eq:expansion_aut}) into eq.
(\ref{eq:invariance_equation_aut}), we equate the terms of equal
polynomial order to obtain a set of \emph{cohomological equations
}\cite{haro16} to be solved recursively in $j$ for the unknown
matrices $\mathbf{W}_{j}^{0},\mathbf{N}_{j}^{0}$. At the $j=1$ step
in this recursion, we have the cohomological equation

\begin{equation}
\mathbf{W}_{1}^{0}\mathbf{N}_{1}^{0}=\boldsymbol{\Lambda}\mathbf{W}_{1}^{0},
\end{equation}
which is solved by 
\begin{equation}
\begin{array}{ccc}
\mathbf{W}_{1}^{0}=\left[\begin{array}{c}
\mathbf{I}\\
\mathbf{0}
\end{array}\right], & \mathbf{N}_{1}^{0}=\mathrm{\boldsymbol{\Lambda}}_{m}=\mathrm{diag}(\boldsymbol{\lambda}_{m}), & \boldsymbol{\lambda}_{m}=(\lambda_{1},\bar{\lambda}_{1},\lambda_{2},\bar{\lambda}_{2},...,\lambda_{m},\bar{\lambda}_{m})^{T}.\end{array}
\end{equation}
Noting that 
\begin{equation}
D\mathbf{w}_{0}(\mathbf{z})\mathbf{N}_{1}^{0}\mathbf{z}={\displaystyle \sum_{j=1}^{M}\mathbf{W}_{j}^{0}(D\mathbf{z}^{j})\mathbf{N}_{1}^{0}\mathbf{z}}={\displaystyle \sum_{j=1}^{M}\mathbf{W}_{j}^{0}\mathrm{\boldsymbol{\Lambda}}_{m}^{j}\mathbf{z}^{j}},
\end{equation}
where $\mathrm{\boldsymbol{\Lambda}}_{m}^{j}=\mathrm{diag}((\mathbf{S}_{2m}^{j})^{T}\boldsymbol{\lambda}_{m})$,
we can rewrite eq. (\ref{eq:invariance_equation_aut}) as 
\begin{equation}
{\displaystyle \sum_{j=1}^{M}\left(\boldsymbol{\Lambda}\mathbf{W}_{j}^{0}-\mathbf{W}_{j}^{0}\mathrm{\boldsymbol{\Lambda}}_{m}^{j}-\mathbf{W}_{1}^{0}\mathbf{N}_{j}^{0}\right)\mathbf{z}^{j}}=\left(\sum_{j=2}^{M}\mathbf{W}_{j}^{0}(D\mathbf{z}^{j})\right)\left(\sum_{j=2}^{M}\mathbf{N}_{j}^{0}\mathbf{z}^{j}\right)-{\displaystyle \sum_{j=2}^{M_{g}}\mathbf{G}_{j}^{0}}\left({\displaystyle \sum_{j=1}^{M}\mathbf{W}_{j}^{0}\mathbf{z}^{j}}\right)^{j}.\label{eq:invariance_equation_aut_nonlin}
\end{equation}

We observe that, in the right hand side of eq. (\ref{eq:invariance_equation_aut_nonlin}),
the coefficients for the monomials of order $k$ only depend on the
coefficients in the matrices $\mathbf{W}_{j}^{0},\mathbf{N}_{j}^{0}$
for $j=1,2,...,k-1$. This motivates us to define 
\begin{equation}
{\displaystyle \sum_{j=2}^{M}\mathbf{B}_{j}^{0}\mathbf{z}^{j}}=\left(\sum_{j=2}^{M}\mathbf{W}_{j}^{0}(D\mathbf{z}^{j})\right)\left(\sum_{j=2}^{M}\mathbf{N}_{j}^{0}\mathbf{z}^{j}\right)-{\displaystyle \sum_{j=2}^{M_{g}}\mathbf{G}_{j}^{0}}\left({\displaystyle \sum_{j=1}^{M}\mathbf{W}_{j}^{0}\mathbf{z}^{j}}\right)^{j},\label{eq:invariance_equation_aut_nonlin_known}
\end{equation}
where, for example, $\mathbf{B}_{2}^{0}\mathbf{z}^{2}=\mathbf{G}_{2}^{0}(\mathbf{W}_{1}^{0}\mathbf{z})^{2}$.
Therefore, the coefficients at order $k$ must satisfy the equation

\begin{equation}
\boldsymbol{\Lambda}\mathbf{W}_{k}^{0}-\mathbf{W}_{k}^{0}\mathrm{\boldsymbol{\Lambda}}_{m}^{k}=\mathbf{W}_{1}^{0}\mathbf{N}_{k}^{0}+\mathbf{B}_{k}^{0}.\label{eq:cohomological}
\end{equation}

Using the notation 
\begin{equation}
\begin{array}{ccc}
\boldsymbol{\Lambda}=\left[\begin{array}{cc}
\boldsymbol{\Lambda}_{m} & \mathbf{0}\\
\mathbf{0} & \boldsymbol{\Lambda}_{out}
\end{array}\right], & \mathbf{W}_{k}^{0}=\left[\begin{array}{c}
\mathbf{W}_{k,in}^{0}\\
\mathbf{W}_{k,out}^{0}
\end{array}\right], & \mathbf{B}_{k}^{0}=\left[\begin{array}{c}
\mathbf{B}_{k,in}^{0}\\
\mathbf{B}_{k,out}^{0}
\end{array}\right]\end{array},
\end{equation}
we decouple the matrix equation (\ref{eq:cohomological}) as 
\begin{equation}
\begin{cases}
\boldsymbol{\Lambda}_{m}\mathbf{W}_{k,in}^{0}-\mathbf{W}_{k,in}^{0}\mathrm{\boldsymbol{\Lambda}}_{m}^{k}=\mathbf{W}_{1}^{0}\mathbf{N}_{k}^{0}+\mathbf{B}_{k,in}^{0},\\
\boldsymbol{\Lambda}_{out}\mathbf{W}_{k,out}^{0}-\mathbf{W}_{k,out}^{0}\mathrm{\boldsymbol{\Lambda}}_{m}^{k}=\mathbf{B}_{k,out}^{0}.
\end{cases}\label{eq:decoupled cohomological equations}
\end{equation}
From the second set of equations in \eqref{eq:decoupled cohomological equations},
we find $\mathbf{W}_{k,out}^{0}$ to be 
\begin{equation}
\left(\mathbf{W}_{k,out}^{0}\right)_{r,s}=\frac{\left(\mathbf{B}_{k,out}^{0}\right)_{r,s}}{\left(\boldsymbol{\Lambda}_{out}\right)_{r,r}-\left(\mathrm{\boldsymbol{\Lambda}}_{m}^{k}\right)_{s,s}},
\end{equation}
because the denominator of this expression is nonzero by the non-resonance
assumption for the existence of autonomous SSMs (see the \emph{Methods}
section of \cite{cenedese21}). The first set of equations in \eqref{eq:decoupled cohomological equations}
is underdetermined, having more unknowns than equations. We can, therefore,
choose different forms for the reduced dynamics, of which we choose
the extended normal-form style parametrization. This avoids small
denominators in the parametrization and hence enhances its domain
of convergence.

Specifically, we follow the extended normal form principle discussed
in \cite{cenedese21} to obtain 
\begin{equation}
\begin{cases}
\left(\mathbf{W}_{k,in}^{0}\right)_{r,s}=0,\,\,\,\,\,\left(\mathbf{N}_{k}^{0}\right)_{r,s}=-\left(\mathbf{B}_{k,in}^{0}\right)_{r,s}, & \mathrm{if\,\,\,\,\,}\left(\boldsymbol{\Delta}^{k}\right)_{r,s}\leq\delta\\
\left(\mathbf{W}_{k,in}^{0}\right)_{r,s}={\displaystyle \frac{\left(\mathbf{B}_{k,in}^{0}\right)_{r,s}}{\left(\boldsymbol{\Lambda}_{m}\right)_{r,r}-\left(\mathrm{\boldsymbol{\Lambda}}_{m}^{k}\right)_{s,s}}},\,\,\,\,\,\left(\mathbf{N}_{k}^{0}\right)_{r,s}=0, & \mathrm{otherwise}.
\end{cases}
\end{equation}
Here $\delta$ is a small tolerance parameter defining near-resonances
among the linearized frequencies related to the SSM, and $\left(\boldsymbol{\Delta}^{k}\right)_{r,s}$
is defined by 
\begin{equation}
\left(\boldsymbol{\Delta}^{k}\right)_{r,s}=\left(\mathrm{Im}\boldsymbol{\Lambda}_{m}\right)_{r,r}-\sum_{l=1}^{2m}\left(\mathrm{Im}\boldsymbol{\Lambda}_{m}\right)_{l,l}\left(\mathbf{S}_{2m}^{k}\right)_{l,s}.
\end{equation}
We perform similar calculations recursively for $j>1$ to obtain the
higher-order polynomial terms in the parametrization and the reduced
dynamics defined in eq. \eqref{eq:expansion_aut}.

\subsection{Solving the non-autonomous invariance equation}

The terms of order $\mathcal{O}(\epsilon)$ in the invariance equation
(\ref{eq:invariance_equation}) must satisfy 
\[
D\mathbf{w}_{0}(\mathbf{z})\mathbf{n}_{1}(\mathbf{z},\boldsymbol{\varphi})+D_{\mathbf{z}}\mathbf{w}_{1}(\mathbf{z},\boldsymbol{\varphi})\mathbf{n}_{0}(\mathbf{z})+D_{\mathbf{\boldsymbol{\varphi}}}\mathbf{w}_{1}(\mathbf{z},\boldsymbol{\varphi})\boldsymbol{\Omega}=\boldsymbol{\Lambda}\mathbf{w}_{1}(\mathbf{z},\boldsymbol{\varphi})+D\mathbf{g}_{0}(\mathbf{w}_{0}(\mathbf{z}))\mathbf{w}_{1}(\mathbf{z},\boldsymbol{\varphi})+\mathbf{g}_{1}(\boldsymbol{\varphi}).
\]
The $\mathcal{O}(1)$ terms of $\mathbf{w}_{1}$ and $\mathbf{n}_{1}$
then satisfy 
\begin{equation}
\boldsymbol{\Lambda}\mathbf{w}_{1}(\mathbf{0},\boldsymbol{\varphi})-D_{\mathbf{\boldsymbol{\varphi}}}\mathbf{w}_{1}(\mathbf{0},\boldsymbol{\varphi})\boldsymbol{\Omega}=\mathbf{W}_{1}^{0}\mathbf{n}_{1}(\mathbf{0},\boldsymbol{\varphi})-\mathbf{g}_{1}(\boldsymbol{\varphi}).\label{eq:invariance_equation_nonaut_zero}
\end{equation}
The last forcing term, 
\begin{equation}
\mathbf{g}_{1}(\boldsymbol{\varphi})={\displaystyle \sum_{|\mathbf{k}|=1}^{\infty}\mathbf{T}^{-1}\mathbf{f}_{\mathbf{k}}^{1}e^{+i\langle\mathbf{k},\boldsymbol{\varphi}\rangle}+}\mathbf{T}^{-1}\bar{\mathbf{f}}_{\mathbf{k}}^{1}e^{-i\langle\mathbf{k},\boldsymbol{\varphi}\rangle}={\displaystyle \sum_{|\mathbf{k}|=1}^{\infty}\mathbf{g}_{\mathbf{k}}^{+}e^{i\langle\mathbf{k},\boldsymbol{\varphi}\rangle}+\mathbf{g}_{\mathbf{k}}^{-}e^{-i\langle\mathbf{k},\boldsymbol{\varphi}\rangle}},\label{eq:forcing_definition_modal}
\end{equation}
satisfies 
\begin{equation}
\begin{array}{ccc}
\left(\bar{\mathbf{g}}_{\mathbf{k}}^{+}\right)_{2j-1}=\left(\mathbf{g}_{\mathbf{k}}^{-}\right)_{2j}, & \left(\bar{\mathbf{g}}_{\mathbf{k}}^{-}\right)_{2j-1}=\left(\mathbf{g}_{\mathbf{k}}^{+}\right)_{2j}, & j=1,2,...,\,m,\end{array}\label{eq:modal_forcing_prop}
\end{equation}
due to the complex conjugate columns of the matrix $\mathbf{T}$.

Substituting the expansions 
\begin{equation}
\begin{array}{c}
{\displaystyle {\displaystyle \mathbf{w}_{1}(\mathbf{z},\boldsymbol{\varphi})=\sum_{|\mathbf{k}|=1}^{\infty}\mathbf{w}_{\mathbf{k}}^{+}e^{i\langle\mathbf{k},\boldsymbol{\varphi}\rangle}+\mathbf{w}_{\mathbf{k}}^{-}e^{-i\langle\mathbf{k},\boldsymbol{\varphi}\rangle}}+\mathcal{O}(\|\mathbf{z}\|),}\\
{\displaystyle \mathbf{n}_{1}(\mathbf{z},\boldsymbol{\varphi})=\sum_{|\mathbf{k}|=1}^{\infty}\mathbf{n}_{\mathbf{k}}^{+}e^{i\langle\mathbf{k},\boldsymbol{\varphi}\rangle}+\mathbf{n}_{\mathbf{k}}^{-}e^{-i\langle\mathbf{k},\boldsymbol{\varphi}\rangle}+\mathcal{O}(\|\mathbf{z}\|),}
\end{array}
\end{equation}
and (\ref{eq:forcing_definition_modal}) into eq. (\ref{eq:invariance_equation_nonaut_zero}),
we solve for $\mathbf{w}_{\mathbf{k}}^{+}$ and $\mathbf{w}_{\mathbf{k}}^{-}$
to obtain 
\begin{equation}
\begin{cases}
{\displaystyle \left(\mathbf{w}_{\mathbf{k}}^{+}\right)_{r}=\frac{\left(\mathbf{n}_{\mathbf{k}}^{+}\right)_{r}-\left(\mathbf{g}_{\mathbf{k}}^{+}\right)_{r}}{\left(\boldsymbol{\Lambda}_{m}\right)_{r,r}-i\langle\mathbf{k},\boldsymbol{\Omega}\rangle},\,\,\,\,\,\,\left(\mathbf{w}_{\mathbf{k}}^{-}\right)_{r}=\frac{\left(\mathbf{n}_{\mathbf{k}}^{-}\right)_{r}-\left(\mathbf{g}_{\mathbf{k}}^{-}\right)_{r}}{\left(\boldsymbol{\Lambda}_{m}\right)_{r,r}+i\langle\mathbf{k},\boldsymbol{\Omega}\rangle},} & \mathrm{if\,\,\,\,\,}r\leq2m\\
{\displaystyle \left(\mathbf{w}_{\mathbf{k}}^{\pm}\right)_{r}=\frac{-\left(\mathbf{g}_{\mathbf{k}}^{\pm}\right)_{r}}{\left(\boldsymbol{\Lambda}_{m}\right)_{r,r}\mp i\langle\mathbf{k},\boldsymbol{\Omega}\rangle}} & \mathrm{otherwise}.
\end{cases}\label{eq:solution_ienaz}
\end{equation}
As previously, the first $2m$ equations in eq. (\ref{eq:solution_ienaz})
are underdetermined. Their simplest meaningful solution yields the
normalized reduced dynamics. This solution is influenced by possible
resonances with the external forcing. Those resonances appear in weakly
damped systems for a certain SSM mode $j\leq m$ whenever $\left(\mathrm{Im}\boldsymbol{\Lambda}_{m}\right)_{2j-1,2j-1}\approx\langle\mathbf{k},\boldsymbol{\Omega}\rangle$
or $\left(\mathrm{Im}\boldsymbol{\Lambda}_{m}\right)_{2j,2j}\approx\langle\mathbf{k},\boldsymbol{\Omega}\rangle$,
where we have $\left(\mathrm{Im}\boldsymbol{\Lambda}_{m}\right)_{2j-1,2j-1}=-\left(\mathrm{Im}\boldsymbol{\Lambda}_{m}\right)_{2j,2j}$.
If there are no such resonances, the forcing vectors in the reduced
dynamics $\mathbf{n}_{\mathbf{k}}^{+},\mathbf{n}_{\mathbf{k}}^{-}$
can be chosen zero. Hence, the linearizable forced state appears in
the parametrization but not in the reduced dynamics, at least not
up to the current order of expansion. In contrast, non-linearizable
forced states may arise due to resonances. As an alternative, we could
retain all forcing terms in the reduced dynamics by setting $\left(\mathbf{n}_{\mathbf{k}}^{-}\right)_{r}=\left(\mathbf{g}_{\mathbf{k}}^{-}\right)_{r}$
and $\left(\mathbf{w}_{\mathbf{k}}^{-}\right)_{r}=0$ for $r\leq2m$,
but this would lead to unnecessary complexity in the reduced model.

Considering a specific mode $j$, we define the sets $K_{j}^{\pm}\subset\mathbb{Z}^{l}$
of resonant indexes as in eq. (\ref{eq:resonant_index}). Adopting
again the normal-form style parametrization, we impose the relationships
\begin{equation}
\begin{cases}
\left(\mathbf{n}_{\mathbf{k}}^{+}\right)_{2j-1}=\left(\mathbf{g}_{\mathbf{k}}^{+}\right)_{2j-1},\,\,\,\left(\mathbf{n}_{\mathbf{k}}^{-}\right)_{2j}=\left(\mathbf{g}_{\mathbf{k}}^{-}\right)_{2j},\,\,\,\left(\mathbf{n}_{\mathbf{k}}^{-}\right)_{2j-1}=\left(\mathbf{n}_{\mathbf{k}}^{+}\right)_{2j}=0, & \mathrm{if}\,\,\,\mathbf{k}\in K_{j}^{+},\\
\left(\mathbf{n}_{\mathbf{k}}^{-}\right)_{2j-1}=\left(\mathbf{g}_{\mathbf{k}}^{-}\right)_{2j-1},\,\,\,\left(\mathbf{n}_{\mathbf{k}}^{+}\right)_{2j}=\left(\mathbf{g}_{\mathbf{k}}^{+}\right)_{2j},\,\,\,\left(\mathbf{n}_{\mathbf{k}}^{+}\right)_{2j-1}=\left(\mathbf{n}_{\mathbf{k}}^{-}\right)_{2j}=0, & \mathrm{if}\,\,\,\mathbf{k}\in K_{j}^{-},
\end{cases}
\end{equation}
while the parametrization terms follow from (\ref{eq:solution_ienaz}).
We note that $K_{j}^{+}\cap K_{j}^{-}=\emptyset$.

\subsection{SSM-reduced dynamics in polar coordinates}

For the $j^{th}$ mode, we find that $\left(\mathbf{n}_{0}(\mathbf{z})\right)_{2j}$
is the complex conjugate of $\left(\mathbf{n}_{0}(\mathbf{z})\right)_{2j-1}$.
Thus, letting $g_{j,\mathbf{k}}^{\pm}=\left(\mathbf{g}_{\mathbf{k}}^{\pm}\right)_{2j-1}$
so that $\left(\mathbf{g}_{\mathbf{k}}^{\mp}\right)_{2j}=\bar{g}_{j,\mathbf{k}}$
as in eq. (\ref{eq:modal_forcing_prop}), we obtain the reduced dynamics
on the SSM in the form 
\begin{equation}
\begin{aligned}\dot{z}_{j} & =\left(\mathbf{n}_{0}(\mathbf{z})\right)_{2j-1}+\sum_{\mathbf{k}\in K_{j}^{+}}\epsilon g_{j,\mathbf{k}}^{+}e^{i\langle\mathbf{k},\boldsymbol{\varphi}\rangle}+\sum_{\mathbf{k}\in K_{j}^{-}}\epsilon g_{j,\mathbf{k}}^{-}e^{-i\langle\mathbf{k},\boldsymbol{\varphi}\rangle}+\mathcal{O}(\epsilon\|\mathbf{z}\|),\\
\dot{\bar{z}}_{j} & =\left(\bar{\mathbf{n}}{}_{0}(\mathbf{z})\right)_{2j-1}+\sum_{\mathbf{k}\in K_{j}^{+}}\epsilon\bar{g}_{j,\mathbf{k}}^{+}e^{-i\langle\mathbf{k},\boldsymbol{\varphi}\rangle}+\sum_{\mathbf{k}\in K_{j}^{-}}\epsilon\bar{g}_{j,\mathbf{k}}^{-}e^{i\langle\mathbf{k},\boldsymbol{\varphi}\rangle}+\mathcal{O}(\epsilon\|\mathbf{z}\|).
\end{aligned}
\label{eq: complex SSM-reduced dynamics}
\end{equation}
Introducing the polar coordinates $z_{j}=\rho_{j}e^{i\theta_{j}}$
and setting $g_{j,\mathbf{k}}^{\pm}=|g_{j,\mathbf{k}}^{\pm}|e^{\pm i(\phi_{j,\mathbf{k}}+\pi/2)}$,
we obtain from \eqref{eq: complex SSM-reduced dynamics} the equations
\begin{equation}
\begin{aligned}\dot{\rho}_{j}+i\rho_{j}\dot{\theta}_{j} & =e^{-i\theta_{j}}\left(\mathbf{n}_{0}(\mathbf{z})\right)_{2j-1}+\sum_{\mathbf{k}\in K_{j}^{\pm}}\epsilon|g_{j,\mathbf{k}}^{\pm}|e^{\pm i(\langle\mathbf{k},\boldsymbol{\varphi}\rangle+\phi_{j,\mathbf{k}}\mp\theta_{j}+\pi/2)}+\mathcal{O}(\epsilon\|\mathbf{z}\|),\\
\dot{\rho}_{j}-i\rho_{j}\dot{\theta}_{j} & =e^{i\theta_{j}}\left(\bar{\mathbf{n}}{}_{0}(\mathbf{z})\right)_{2j-1}+\sum_{\mathbf{k}\in K_{j}^{\pm}}\epsilon|g_{j,\mathbf{k}}^{+}|e^{\mp i(\langle\mathbf{k},\boldsymbol{\varphi}\rangle+\phi_{j,\mathbf{k}}\mp\theta_{j}+\pi/2)}+\mathcal{O}(\epsilon\|\mathbf{z}\|),
\end{aligned}
\label{eq:first polar form of SSM-reduced dynamics}
\end{equation}
where we have grouped together the $\pm$ sets for notational ease,
i.e., 
\[
\begin{array}{cc}
{\displaystyle \sum_{\mathbf{k}\in K_{j}^{\pm}}\epsilon|g_{j,\mathbf{k}}^{\pm}|e^{\pm i(\langle\mathbf{k},\boldsymbol{\varphi}\rangle+\phi_{j,\mathbf{k}}\mp\theta_{j}+\pi/2)}=} & {\displaystyle \sum_{\mathbf{k}\in K_{j}^{+}}\epsilon|g_{j,\mathbf{k}}^{+}|e^{i(\langle\mathbf{k},\boldsymbol{\varphi}\rangle+\phi_{j,\mathbf{k}}-\theta_{j}+\pi/2)}+}\\
 & {\displaystyle +\sum_{\mathbf{k}\in K_{j}^{-}}\epsilon|g_{j,\mathbf{k}}^{-}|e^{-i(\langle\mathbf{k},\boldsymbol{\varphi}\rangle+\phi_{j,\mathbf{k}}+\theta_{j}+\pi/2)}.}
\end{array}
\]

By separating the time derivatives of the amplitude and phase variables
in \eqref{eq:first polar form of SSM-reduced dynamics}, we obtain
the normal form 
\begin{equation}
\begin{aligned}\dot{\rho}_{j} & =\mathrm{Re}\left(e^{-i\theta_{j}}\left(\mathbf{n}_{0}(\mathbf{z})\right)_{2j-1}\right)-\sum_{\mathbf{k}\in K_{j}^{\pm}}\epsilon|g_{j,\mathbf{k}}^{\pm}|\sin(\langle\mathbf{k},\boldsymbol{\varphi}\rangle+\phi_{j,\mathbf{k}}\mp\theta_{j})+\mathcal{O}(\epsilon\|\boldsymbol{\rho}\|),\\
\rho_{j}\dot{\theta}_{j} & =\mathrm{Im}\left(e^{-i\theta_{j}}\left(\mathbf{n}_{0}(\mathbf{z})\right)_{2j-1}\right)+\sum_{\mathbf{k}\in K_{j}^{\pm}}\epsilon|g_{j,\mathbf{k}}^{\pm}|\cos(\langle\mathbf{k},\boldsymbol{\varphi}\rangle+\phi_{j,\mathbf{k}}\mp\theta_{j})+\mathcal{O}(\epsilon\|\boldsymbol{\rho}\|).
\end{aligned}
\label{eq:second polar form of SSM-reduced dynamics}
\end{equation}
With the definitions 
\begin{equation}
\begin{aligned}\alpha_{j}(\boldsymbol{\rho},\boldsymbol{\theta}) & =\mathrm{Re}\left(\frac{\left(\mathbf{n}_{0}(\mathbf{z})\right)_{2j-1}}{z_{j}}\right),\\
\omega_{j}(\boldsymbol{\rho},\boldsymbol{\theta}) & =\mathrm{Im}\left(\frac{\left(\mathbf{n}_{0}(\mathbf{z})\right)_{2j-1}}{z_{j}}\right),
\end{aligned}
\label{eq:amp_phase_dyn}
\end{equation}
and with the rescaling $f_{j,\mathbf{k}}=\epsilon|g_{j,\mathbf{k}}^{\pm}|$,
eq. \eqref{eq:second polar form of SSM-reduced dynamics} provides
the detailed form of the polar normal form given in the methods section
of \cite{cenedese21}. We also note that 
\[
{\displaystyle \lim_{\|\boldsymbol{\rho}\|\rightarrow0}\alpha_{j}(\boldsymbol{\rho},\boldsymbol{\theta})+i\omega_{j}(\boldsymbol{\rho},\boldsymbol{\theta})}=\left(\boldsymbol{\Lambda}_{m}\right)_{2j-1,2j-1}.
\]
For a two-dimensional SSM, periodically forced ($l=1$) with a single
sinusoidal term near resonance, we have $K^{+}=\left\{ 1\right\} $,
$K^{-}=\emptyset$. The amplitude and phase dynamics in (\ref{eq:amp_phase_dyn})
only depends on $\rho$, so that we recover eq. (7) of \cite{cenedese21}
by setting $\theta=\Omega t+\phi_{1,1}+\psi$.

The model reduction we have outlined for periodic and quasiperiodic
SSMs features an autonomous core plus the leading-order forcing term,
providing an overall $\mathcal{O}(\epsilon\|\mathbf{\boldsymbol{\rho}}\|)$
accuracy. More accurate approximations can also be derived including
the remaining $\mathcal{O}(\epsilon)$ terms that also depend on higher
powers of the amplitudes $\boldsymbol{\rho}$ or higher powers of
$\epsilon$. While the leading-order reduced dynamics, $\mathbf{n}_{1}(\boldsymbol{\varphi})$,
can already identify coexisting isolated steady states, the parametrization
correction term, $\mathbf{w}_{1}(\boldsymbol{\varphi})$, reveals
an additional, small quasiperiodic modulation to these autonomous
steady states. These two terms both depend on the forcing coefficients
$g_{j,\mathbf{k}}$, which can be identified via calibration to experiments.
If information on modes outside those related to the slow $2m$-dimensional
SSM is not available, then the reduced-order model is still accurate
up to a small quasiperiodic correction if we set $\left(\mathbf{w}_{\mathbf{k}}^{\pm}\right)_{2j-1}=\left(\mathbf{w}_{\mathbf{k}}^{\pm}\right)_{2j}=0$
for $j>m$. However, if there are resonances involving those additional
modes, it is advisable to increase the dimension of the SSM by including
them. This will generally require the collection of further data with
additional modal content.
